\documentclass[11pt]{article}
\pdfoutput=1 
\usepackage{amssymb}
\usepackage[margin=1in]{geometry}

\usepackage{lineno}
\usepackage{graphicx}

\usepackage{moreverb,url}
\usepackage{makecell}
\usepackage[breaklinks=true, colorlinks=true, linkcolor={blue!90!black}, citecolor={blue!90!black}, urlcolor={blue!90!black}]{hyperref}

\newcommand\BibTeX{{\rmfamily B\kern-.05em \textsc{i\kern-.025em b}\kern-.08em
T\kern-.1667em\lower.7ex\hbox{E}\kern-.125emX}}

\newcommand{\defrefstwo}[2]{{\rm \textsc{Definitions}~\labelcref{#1}}~and~{\rm \labelcref{#2}}}

\usepackage[ruled,vlined]{algorithm2e}
\usepackage{algpseudocode}
\usepackage{pifont}
\usepackage{array}
\usepackage{tabularx}
\usepackage{mathptmx}

\DeclareMathAlphabet{\mathcal}{OMS}{cmsy}{m}{n}

\newcommand{\cmark}{\ding{51}}%
\newcommand{\xmark}{\ding{55}}%

\SetCommentSty{mycommfont}

\usepackage{geometry,xcolor,tabularx,nccmath}
\usepackage{amsmath,amsthm,amssymb,enumitem}
\usepackage{booktabs,caption,subcaption,mathtools,multirow}
\usepackage{graphicx,tabulary}
\usepackage{rotating}

\usepackage[utf8]{inputenc}
\usepackage{amsfonts}
\usepackage{fullpage}

\newcommand{\sC}{\mathcal{C}}

\newcommand{\sI}{\mathcal{I}}
\newcommand{\sJ}{\mathcal{J}}
\newcommand{\sK}{\mathcal{K}}

\newcommand{\sP}{\mathcal{P}}

\newcommand{\sS}{\mathcal{S}}

\newcommand{\sX}{\mathcal{X}}

\newcommand{\bigM}{\mathbb{M}}
\newcommand{\norm}[1]{\lvert #1 \rvert}

\newcommand{\f}{f_{u^*,v^*,i^*,j^*}}

\usepackage[colorinlistoftodos]{todonotes}

\usepackage{etoolbox}

\usepackage[nameinlink,capitalise]{cleveref}
\hypersetup{colorlinks={true},allcolors={blue}}

\crefformat{section}{\S#2#1#3}
\crefformat{subsection}{\S#2#1#3}
\crefformat{equation}{#2Equation~\textup{(#1)}#3}
\crefformat{cons}{#2Constraint~\textup{(#1)}#3}
\labelcrefformat{cons}{\textup{#2(#1)#3}}
\crefformat{consset}{#2Constraints~\textup{(#1)}#3} 

\labelcrefformat{consset}{\textup{#2(#1)#3}}
\crefformat{func}{#2Function~\textup{(#1)}#3}
\labelcrefformat{func}{#2Function~\textup{(#1)}#3}
\crefformat{algo}{#2Algorithm~\textup{#1}#3}
\newcommand{\algorefstwo}[2]{{\rm Algorithms~\labelcref{#1}}~and~{\rm \labelcref{#2}}}
\crefformat{objfunc}{#2Objective function~\textup{(#1)}#3}
\labelcrefformat{objfunc}{\textup{#2(#1)#3}}
\crefformat{ch}{#2Chapter~\textup{#1}#3}
\crefformat{sec}{#2Section~\textup{#1}#3}
\crefformat{tab}{#2Table~\textup{#1}#3}
\crefformat{fig}{#2Figure~\textup{#1}#3}
\crefdefaultlabelformat{#2#1#3} 
\Crefname{figure}{Figure}{Figures}
\crefformat{eq}{#2equation~\textup{(#1)}#3}
\labelcrefformat{eq}{\textup{(#2#1#3)}}
\crefformat{theo}{#2Theorem~\textup{#1}#3}
\crefformat{lem}{#2Lemma~\textup{#1}#3}
\crefformat{def}{#2Definition~\textup{#1}#3}
\crefformat{app}{#2\textup{\!#1}#3}
\crefformat{proper}{#2\textsc{Property}~\textup{\textsc #1}#3}
\newcommand{\proprefstwo}[2]{{\rm \textsc{Properties}~\labelcref{#1}}~and~{\rm \labelcref{#2}}}
\crefformat{item}{#2Step~\textup{#1}#3}
\crefformat{scenario}{#2Scenario~\textup{#1}#3}
 \AtBeginDocument{%
    \crefname{figure}{Figure}{figures}%
}
\let\origref\cref
\def\cref#1{\origref{#1}}

\expandafter\def\expandafter\UrlBreaks\expandafter{\UrlBreaks
  \do\a\do\b\do\c\do\d\do\e\do\f\do\g\do\h\do\i\do\j%
  \do\k\do\l\do\m\do\n\do\o\do\p\do\q\do\r\do\s\do\t%
  \do\u\do\v\do\w\do\x\do\y\do\z\do\A\do\B\do\C\do\D%
  \do\E\do\F\do\G\do\H\do\I\do\J\do\K\do\L\do\M\do\N%
  \do\O\do\P\do\Q\do\R\do\S\do\T\do\U\do\V\do\W\do\X%
  \do\Y\do\Z}

\usepackage{amsmath,cases,eqparbox,xparse}

\makeatletter
\NewDocumentCommand{\eqmathbox}{o O{c} m}{%
  \IfValueTF{#1}
    {\def\eqmathbox@##1##2{\eqmakebox[#1][#2]{$##1##2$}}}
    {\def\eqmathbox@##1##2{\eqmakebox{$##1##2$}}}
  \mathpalette\eqmathbox@{#3}
}
\makeatother

\makeatletter
\let\@msm@th@eqref\eqref
\renewcommand{\eqref}[1]{%
  \begingroup
  \leavevmode
  \color{blue}%
  \hypersetup{linkbordercolor=[named]{blue}}%
  \@msm@th@eqref{#1}%
  \endgroup
}
\makeatother

\usepackage[normalem]{ulem} 
 
\crefformat{item}{#2Step~\textup{#1}#3}

\crefname{defi}{definition}{definitions}
\Crefname{defi}{Definition}{Definitions}
\crefname{lemma}{lemma}{lemmas}
\Crefname{lemma}{Lemma}{Lemmas}
\crefname{assumption}{assumption}{assumptions}
\Crefname{assumption}{Assumption}{Assumptions}

\usepackage{multicol}

\usepackage{authblk}
\providecommand{\keywords}[1]{\noindent\textbf{Keywords:} #1}

\usepackage[authoryear]{natbib}

\newtheorem{property}{Property}

\crefformat{proper}{#2\textsc{Property}~\textup{\textsc #1}#3}

\newtheorem{definition}{Definition}
\crefformat{def}{#2\textsc{Definition}~\textup{\textsc #1}#3}

\begin{document}

\title{Problem of Locating and Allocating Charging Equipment for Battery Electric Buses under Stochastic Charging Demand}


\author[1]{Sadjad Bazarnovi}
\author[2,*]{Taner Cokyasar}
\author[2]{Omer Verbas}
\author[1]{Abolfazl (Kouros) Mohammadian}
\affil[1]{University of Illinois Chicago, 842 W Taylor St, Chicago IL 60607, USA }
\affil[2]{Argonne National Laboratory, 9700 S. Cass Ave., Lemont, IL 60439, USA }
\affil[*]{Corresponding Author}

\maketitle


\begin{abstract}
Bus electrification plays a crucial role in advancing urban transportation sustainability. Battery Electric Buses (BEBs), however, often need recharging, making the Problem of Locating and Allocating Charging Equipment for BEBs (PLACE-BEB) essential for efficient operations. This study proposes an optimization framework to solve the PLACE-BEB by determining the optimal placement of charger types at potential locations under the stochastic charging demand. Leveraging the existing stochastic location literature, we develop a Mixed-Integer Non-Linear Program (MINLP) to model the problem. To solve this problem, we develop an exact solution method that minimizes the costs related to building charging stations, charger allocation, travel to stations, and average queueing and charging times. Queueing dynamics are modeled using an \textit{M}/\textit{M}/\textit{s} queue, with the number of servers at each location treated as a decision variable. To improve scalability, we implement a Simulated Annealing (SA) and a Genetic Algorithm (GA) allowing for efficient solutions to large-scale problems. The computational performance of the methods was thoroughly evaluated, revealing that SA was effective for small-scale problems, while GA outperformed others for large-scale instances. A case study comparing garage-only, other-only, and mixed scenarios, along with joint deployment, highlighted the cost benefits of a collaborative and a comprehensive approach. Sensitivity analyses showed that the waiting time is a key factor to consider in the decision-making.
\end{abstract}


\section{Introduction} \label[sec]{intro}
With more than two-thirds of the world's population predicted to live in cities by 2050 \citep{PERUMAL2022395}, urban transportation networks will play a vital role in enhancing the quality of life and the functionality of cities by providing essential mobility and accessibility to residents. Public transportation, in particular, is key to reducing traffic congestion, lowering emissions, and providing equitable access to transportation. However, transportation substantially contributes to air pollution, with the transportation sector being one of the largest producers of carbon dioxide \citep{TZAMAKOS2023291, USLU2021102645}. Reducing emissions is critical in the fight against climate change, and one promising solution is the adoption of vehicles powered by alternative fuels, such as electric, hybrid, and hydrogen-powered vehicles \citep{HSU2021103053, DMcCabe2023}.

One significant step towards reducing carbon emissions and noise pollution in urban areas is the electrification of public transportation bus fleets \citep{TZAMAKOS2023291, KGkiotsalitis2023}. Numerous cities worldwide are transitioning to Battery Electric Buses (BEBs) to achieve this goal \citep{TZAMAKOS2023291, DMcCabe2023, HaoHu2022}. For instance, the Chicago Transit Authority (CTA) has committed to transitioning to a fully electric bus fleet by 2040 \citep{authority2022charging}. The benefits of BEBs include zero tailpipe emissions, reduced maintenance costs due to fewer moving parts, decreased noise pollution, reduced vibration, enhanced passenger comfort, and lower fuel and maintenance costs \citep{HaoHu2022, TZAMAKOS2023291, USLU2021102645, DMcCabe2023}. However, BEBs also come with challenges such as limited travel range, long charging times, increased planning complexity, and high initial costs for both the buses and the required charging infrastructure \citep{YiHe2023, KGkiotsalitis2023, TZAMAKOS2023291, DMcCabe2023}.

There are two prevailing approaches to charging BEBs: on-route charging at terminal points ---including a bus garage and other locations--- using fast chargers and overnight charging with slow chargers at the garage. On-route charging is faster but requires more chargers and infrastructure \citep{HaoHu2022}, while slow chargers are cheaper but necessitates a larger fleet to cover all scheduled trips and put significant demand on the power grid at night \citep{PVende2023}. A balanced approach is essential, as relying solely on garage charging would require a larger fleet and cause grid strain, while exclusive reliance on on-route charging would lead to increased costs due to the need for more chargers and stations.

The Problem of Locating and Allocating Charging Equipment for Battery Electric Buses (PLACE-BEB) involves balancing costs of building stations and allocating chargers with operational costs associated with bus waiting times and travel times for charging. Increasing the number of stations can reduce travel costs, as buses can quickly reach nearby stations when they need to charge. Similarly, increasing the number of chargers at each station potentially reduces waiting times, as buses spend less time queueing to charge. On the other hand, adding the possible redundancy comes at a high deployment cost. To this end, a balanced objective approach \citep{bermanandkrass} will equilibrate the level of deployment and the cost of potential queueing and travel.

\begin{figure}[!ht]
    \centering
    \includegraphics[width=1\linewidth]{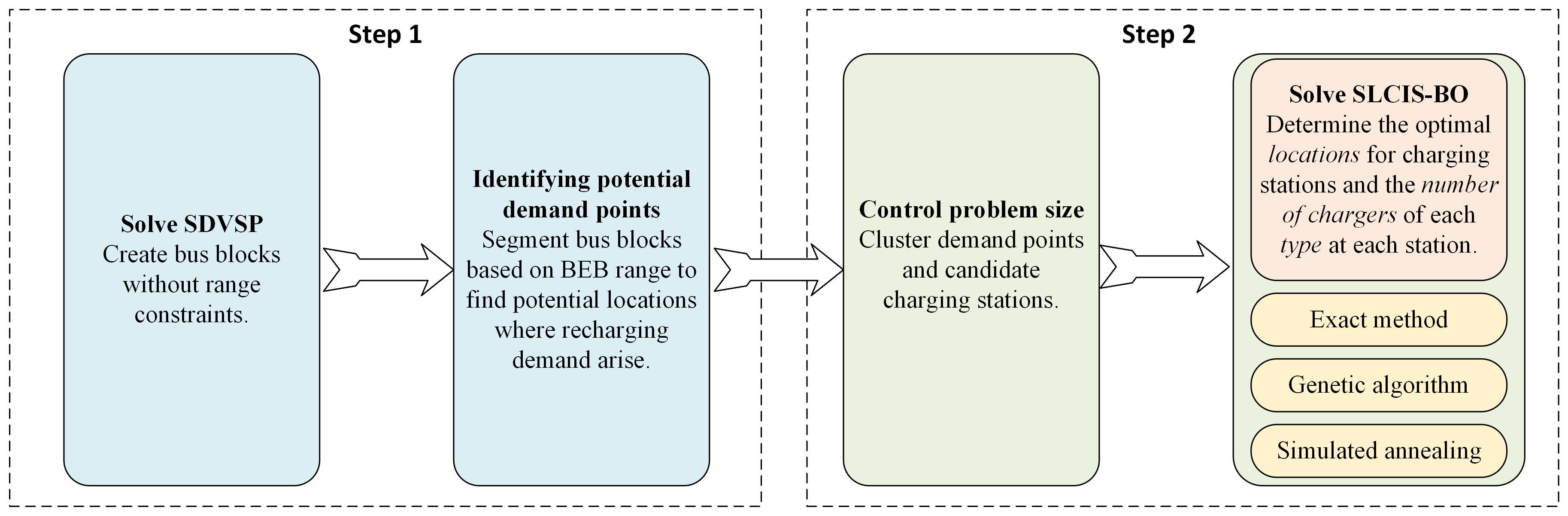}
    \caption{Procedure of the study.}
    \label[fig]{procedure}
\end{figure}

To address the PLACE-BEB, we develop an optimization framework that jointly finds the optimal locations for charging stations and the number and the type of chargers at each candidate station. The framework is composed of two steps, and \cref{procedure} depicts a flowchart of the framework. 

In the first step, we solve the Single Depot Vehicle Scheduling Problem (SDVSP) for buses assuming they are not constrained by an electric driving range. In the SDVSP, revenue-generating timetabled trips are chained into bus schedule \emph{blocks} while minimizing the total intertrip layover time, deadheading time, and the number of blocks. In this step, we identify potential locations where recharging demand emerge by segmenting blocks based on a BEB range and mapping the corresponding trip terminals. 

In the second step, we extend a model belonging to the Stochastic Location Models with Congestion, Immobile Server, and Balanced Objective (SLCIS-BO) literature to find the optimal location, number, and type of chargers \citep{bermanandkrass}. The PLACE-BEB is hence modeled with a Mixed-Integer Non-Linear Program (MINLP). The stochastic nature of the problem is addressed utilizing the queueing theory, and the impact of average waiting times at chargers is incorporated into the location and allocation decision-making. Since the problem scale dealt with is large, this step also involves a post-processing of the solutions obtained through the first step to control the problem size. In the post-processing, recharging demand locations are clustered to control the problem size. The main focus of this paper however is the second step and the three solution methods presented: exact solution (called MINLP), Simulated Annealing (SA) algorithm, and Genetic Algorithm (GA). The first step prepares the input data for the MINLP, and the clustering approach controls the problem size to achieve high quality solutions with the exact solution method under reasonable computational restraints.

The exact solution algorithm is derived through adopting cutting-plane constraints. While the algorithm finds optimal solutions to small-scale problems quickly, large-scale problems require unfavorably high computational times. Therefore, inspired by the existing SLCIS-BO literature, we propose two metaheuristic approaches: an SA and a GA. Comparisons of these methods with exact solutions for smaller instances demonstrate their efficiency and ability to quickly provide good solutions. In our case studies, we focus the public transportation system in the Chicago metropolitan region, specifically CTA and Pace Suburban Bus (henceforth called Pace).

Contributions of the paper are outlined as follows:
\begin{itemize}[noitemsep]
    \item The recharging demand stochasticity is integrated into charger location and allocation decisions. The stochasticity aspect was considered only in a few studies, and this study extends upon the SLCIS-BO literature with an application into the BEB context and a removal of a constraint that allows finding better solutions. We name the constraint as \emph{proximity assignment}, and it ensures assigning each demand point to the closest available supply location. The constraint is further elaborated in \cref{model}.
    \item The literature does not consider a joint optimization of charger location and the number of different charger types in the urban transit electrification context under stochasticity. While integrating both, this study also takes the waiting time into account and fills the gap in the literature. 
    \item Coverage constraints of the classical SLCIS-BO are replaced with subset representations to improve the robustness, and variable charger types are incorporated into the model.
    \item The exact solution method provides a methodological contribution to the literature by being able to address the SLCIS-BO models with or without the proximity assignment constraints.
    \item Meta-heuristic solution approaches are designed to address the scalability concern, making the modeling applicable to real-world problem instances.
    \item Case studies are conducted using real-world bus networks of two agencies operating in the Chicago metropolitan region: CTA and Pace.
\end{itemize}

The remainder of this paper is organized as follows: \cref{lit_rev} reviews relevant literature on optimizing charging station locations and charger allocation. \cref{methodology} formulates the problem. \cref{approaches} presents the solution methods, including an exact method (\cref{exact_method}) and two metaheuristic approaches (\cref{GA} and \cref{SA}). In  \cref{num_exp}, we provide details on our data, present results of a thorough computational experiments by highlighting the capabilities of the solution methods developed, provide managerial insights through case studies comparing various deployment schemes, and pointing out key parameters in decision-making via extensive sensitivity analyses. \cref{conclusion} summarizes our findings and contributions. Finally, figures, tables, and algorithms that are not essential for readability are provided in the \hyperref[figures]{Supplementary Materials}.

\section{Literature review} \label[sec]{lit_rev}
With advancements in Electric Vehicle (EV) technologies and growing concerns about greenhouse gas emissions and the sustainability of transportation networks, many individuals and transportation companies are transitioning from internal combustion engine vehicles to EVs \citep{HaoHu2022, HSU2021103053}. Numerous studies have explored the deployment of charging infrastructure for both private and heavy-duty EVs. The urban transit sector offers significant potential for EV deployment, as buses typically operate on shorter trip distances and fixed routes, making them well-suited to the current limitations of battery range \citep{HSU2021103053}. While many studies have focused on optimizing charging infrastructure for EVs in general \citep{sun2020integrated, xu2020optimal, he2018optimal}, the unique characteristics of BEBs require special attention \citep{HaoHu2022}. This has lead to several studies specifically targeting the placement of BEB charging stations in urban areas, aiming to provide efficient and cost-effective charging solutions while improving service levels for passengers \citep{TZAMAKOS2023291, rogge2018electric, wei2018optimizing, HSU2021103053, HaoHu2022, DMcCabe2023, YiHe2023, an2020battery}.

\cite{he2018optimal} tackled the problem of locating charging stations for EVs through a bi-level modeling approach. Their study addressed the relationship between charging station deployment and route selection. Their findings emphasized that the battery range of EVs and the distances driven significantly influence the optimal placement of charging infrastructure.

\cite{rogge2018electric} developed a model that simultaneously optimizes bus fleet size, vehicle scheduling, charger numbers, and charging schedules for transit systems employing BEBs. They employed a GA to solve this model and applied it in a case study involving 200 trips across three bus routes. The study's limitation to garage charging assumed that BEBs return to the garage for recharging, and restricted its applicability to scenarios with on-route layover charging.

\cite{wei2018optimizing} introduced a spatio-temporal optimization model aimed at minimizing costs related to vehicle procurement and charging station allocation, while adhering to existing bus routes and schedules. Their MILP model identified optimal route assignments, as well as the locations and sizes of on-route layover charging stations and overnight charging sites at garages. They applied this model to the transit network operated by the Utah Transit Authority. However, their assumption that buses can fully charge during any layover period of 10 minutes or more makes their model impractical unless buses have small batteries or use extremely high-power chargers.

\cite{an2020battery} developed a stochastic integer model to solve the problem of optimizing the placement of plug-in chargers and fleet size for BEBs. The study considered time-varying electricity prices and the fluctuations in charge demand due to weather or traffic conditions. The model assumed vehicles may charge only between trips while maintaining the existing diesel bus operation schedule, aggregating demand at terminals based on battery levels. It established discrete one-hour charging time blocks, ensuring the total number of available blocks is not exceeded.

\cite{USLU2021102645} proposed a mixed integer-linear programming model to optimize the placement and capacity of BEB charging stations. The model took into account routes, demand, and driving ranges of BEBs, and incorporates constraints on waiting times using queue theory. Implemented in a case study for intercity bus networks in Türkiye using real-world data, the results highlighted optimal station locations and capacities with minimal costs. However, since their model was designed for a single-line intercity bus service, allowing charging stations to be placed anywhere along the route, urban transit agencies cannot utilize it readily.

\cite{HSU2021103053} solved the problem of locating charging facilities for electrifying urban bus services, highlighting the transition from diesel to electric buses. Their optimization model addressed practical concerns such as fleet size, land acquisition, bus allocation, and deadhead mileage. They proposed a decomposition-based heuristic algorithm for computational efficiency and applied it to a bus operator in Taiwan.

\cite{HaoHu2022} introduced a joint optimization model aimed at placing fast chargers at bus stops and establishing effective charging schedules. Initially, they approached the location problem deterministically, taking into account bus battery capacity and charging requirements. Subsequently, they incorporated uncertainties related to passenger boarding and bus travel duration. Their study applied this model to several bus routes in Sydney, revealing that optimal charging strategies involved frequent charging at intermediate and final stops, leveraging passenger boarding and alighting times. 

\cite{TZAMAKOS2023291} developed a model for optimally locating fast wireless chargers in electric bus networks, accounting for bus delays due to queuing at charging locations. The integer linear programming model minimized investment costs for opportunity charging facilities, employing an \textit{M}/\textit{M}/1 queuing model to incorporate queuing delays. The results indicate the importance of considering waiting times when defining the number of charger stations. However, they only consider opportunity charging during bus service trips and do not include charging during layover times at garages or terminal stations.

\cite{YiHe2023} addressed the combined challenges of charging infrastructure planning, vehicle scheduling, and charging management for BEBs. Their MINLP formulation aimed to minimize the total cost of ownership, employing a genetic algorithm-based approach for solution. The study, applied to a sub-transit network in Salt Lake City, Utah, compared alternative scenarios to the optimal results, demonstrating the model's effectiveness in developing cost-efficient planning strategies.

\cite{DMcCabe2023} presented the BEB optimal charger location model, an MILP approach to optimize the locations and sizes of layover charging stations for BEBs. The model balances infrastructure costs with operational performance, ensuring no bus waits at a charging location. Additionally, the BEB block revision problem model revises vehicle schedules to dispatch backup buses for trips at risk of running out of battery.

The reviewed studies collectively highlight the complexity and importance of strategic planning for charging infrastructure. Key considerations include minimizing costs, ensuring operational efficiency, and accommodating practical constraints such as waiting times, and driving ranges. The models and methodologies developed in these studies provide valuable insights and tools for transit agencies and planners aiming to optimize the deployment of BEB charging stations. \cref{lit_sum} compares this paper with relevant studies in the literature. The table highlights a significant gap: no existing methods optimize both the location of charging stations and the number of different types of chargers at each station while accounting for waiting times in an urban transit network. This paper aims to address this gap by developing an optimization model for the strategic placement of charging stations. Furthermore, the proposed meta-heuristic approaches ensure that the model can be effectively applied to real-world networks.

\section{Methodology} \label[sec]{methodology}
In this section, we first provide a detailed explanation of the problem, followed by an outline of our assumptions. Next, we introduce the notations used in the model. Finally, we present the optimization model, explaining the objective function and all associated constraints.

\subsection{Problem overview} \label{sec:problem_overview}
We begin with reviewing the first step of the PLACE-BEB that prepares inputs to the second step. The classical SDVSP involves creating bus blocks with a given set of timetabled service trips. Each trip represents a service plan with essential spatio-temporal information that are origin, destination, start time, and end time. The origin and destination are also called terminal points. The charger location problem is ideally solved in conjunction with the scheduling problem. However, the electrification occurs over time, that is transit agencies alter their fleets with the electric ones gradually. To this end, a typical BEB deployment behavior is to electrify existing routes and bus schedules rather than re-planning the entire system. To mimic this behavior, we begin with solving SDVSPs for each garage of a given transit agency using their trip information. Hence, the SDVSP solutions yield fossil fuel powered bus blocks. Blocks represent a sequence of consecutive trips between two garage visits \citep{ADavatgari2024}. Using these blocks and a given BEB range, we then partition blocks and find potential demand points from where BEBs could head to charging stations. 

\begin{figure}[!htbp]
    \centering
    \includegraphics[width=\linewidth]{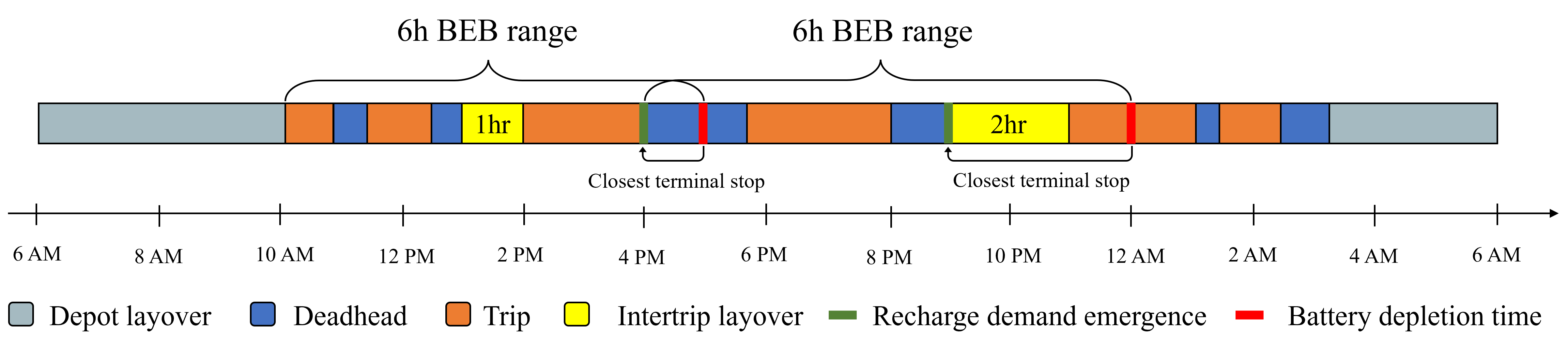}
    \caption{Identifying potential demand points.}
    \label[fig]{blocksegmenting}
\end{figure}

\cref{blocksegmenting} illustrates an example bus block including layover at the garage, deadheading from/to garage and trips, trip timing, layover time between trips, timing of recharge demand, and battery depletion time. Here, the bus starts the service at 6 AM, deadheads to another trip's origin, serves the second trip, deadheads to the third trip's origin, lays over at this terminal for an hour, and serves its third trip. We assume buses do not consume energy during layover. If this block is run by a BEB and assuming a BEB range of six hours, it cannot complete the fourth trip. The red vertical lines in the figure show the battery depletion time based on the BEB range and are called \emph{dividers}. Buses are assumed to charge only when not in service, that is before or after a service trip similar to layover charging presented in \cite{DMcCabe2023}. Thus, although the battery depletes at 5 PM, the charging event will occur at 4 PM at the destination terminal of the third trip, which is the closest terminal stop before the battery depletion time. Assuming the bus fully recharges, the next charging demand arises at around 12 AM.  Conducting a similar logic-based strategy for all bus blocks, we find demand points and the number of charging need at each point.

\cref{DemandFlowChart} depicts the flowchart for positioning the dividers between trips of a BEB. The State-of-Charge (SOC) is a key element in this positioning. It is worthwhile to note that the process in this figure is only one way to estimate potential demand points, and a different method for this estimation could be used. We assume that demand points and the number of charging need at each point are present a priori.

We assume a set of candidate charging stations is given. The set can include any bus stop and garage locations in the network. The goal and the main focus of this study is to determine which stations to deploy a number of chargers with various types, that is the second step of the PLACE-BEB. 

\subsection{Assumptions}\label{sec:assumption}
We develop the model and solution framework based on several key assumptions, drawing on studies that have utilized similar frameworks. It's important to note that some of these assumptions can be adjusted, often requiring minimal or no changes to the model itself. We present these assumptions to inform the reader about the essential data and conditions needed for the model's effectiveness.

\begin{itemize}[noitemsep]
    \item Bus schedules and blocks are provided in advance \citep{DMcCabe2023, HaoHu2022}, which includes detailed information, such as trip start and end times, deadhead and layover times, terminal stops, and garages \citep{DMcCabe2023, HSU2021103053, TZAMAKOS2023291}. If charging demand information through another source is available, these information are not necessary.
    \item BEBs are homogeneous, including battery capacity, range, and efficiency \citep{TZAMAKOS2023291, HSU2021103053}.
    \item BEBs operate within a predefined bound of SOC levels. We assume a range of 10-80\% since energy consumption and gain are linearly correlated to time within this bound \citep{DMcCabe2023, TZAMAKOS2023291}.
    \item Charging demand follows a Poisson process, and charging times are exponentially distributed.
    \item Each vehicle is assigned to a single charging station and single charger type.
    \item Candidate charging stations possess adequate space and power grid capacity to accommodate large number of chargers.
    \item Charging stations can be equipped with multiple charger types \citep{TZAMAKOS2023291}.
\end{itemize}

\subsection{Notations} \label{sec:notation}
We use calligraphic letters to represent sets (e.g., $\sI$), lowercase Roman letters for variables and indices (e.g., $x_{ijk}$), Greek letters for superscripts modifying parameters and variables (e.g., $\delta$), and uppercase Roman letters for parameters (e.g., $T_{ij}$). Blackboard bold letters denote functions, domains, and model names (e.g., $\mathbb{W}$), while blackboard bold letters with a center dot (e.g., $\mathbb{W}(\cdot)$) abbreviate the inputs of the corresponding function. The superscript $*$ (asterisk) indicates optimality, and bold letters denote solution vectors (e.g., $\mathbf{x}$). Notations specific to queueing theory, such as $W_{jk}$ for variables and $\lambda_i$ and $\mu_k$ for parameters, as well as $\epsilon$ for an infinitesimal number and $\bigM$ for a sufficiently large number, are excluded from this notation convention.

\subsection{Mathematical model} \label[sec]{model}
Let $i\in\sI$ denote a set of demand points from where BEBs head to chargers. Candidate charging stations are denoted by $j\in\sJ$, and the objective is to activate a subset of them for charger deployment. The binary variable $y_j=1$ indicates that station $j$ is activated and incurs $C^\phi_j$ cost per time unit. A station from the subset $\sJ_i \subseteq \sJ$ can be visited to charge buses. This subset replaces the commonly adopted \emph{coverage constraints} \citet[Ch. 17, pp. 486]{bermanandkrass} such that additional constraints are not needed in the model to restrict the stations that could be visited. Considering all stations as potential locations to visit for charging both explodes the problem size and is unrealistic. The SOC may not allow visiting any given station in the network. This subset could be refined based on a maximum travel distance/time from demand points to charging stations. We also define $\sI_j$ as the subset of demand points from where BEBs can head to charging at station $j\in\sJ$. \cref{notations} lists all sets, parameters, and variables used in the mathematical model. Unless otherwise noted, all parameters are non-negative real numbers ($\mathbb{R}_{\geq 0}$).

There is a set of $\sK$ charger types with various power outputs to be allocated to stations. The service rate for chargers, i.e., the number of vehicles that can be charged per time unit, is denoted by $\mu_k$ and depends on the initial SOC of vehicles before charging ($Q^\iota$), the targeted SOC level after charging ($Q^\phi$), the battery capacity of vehicles ($B$), and the charger power ($P_k$). Yet, the SOC is not tracked in the stochastic environment, therefore, to estimate \(\mu_k\), all BEBs are assumed to arrive at the charging station with a fixed \(Q^\iota\) and depart with a fixed \(Q^\phi\). Taking the difference between two SOC levels and using $B$ and $P_k$, we can then estimate an expected time to recharge BEBs via charger type $k\in\sK$, $T_k^\rho$. Let binary variable $x_{ijk}=1$ denote that charger type $k\in\sK$ at station $j\in\sJ_i$ is visited to charge a BEB traveling from demand point $i\in\sI$. The number of type $k$ chargers, each incurring $C^\xi_k$ cost per time unit, allocated to station $j$ is denoted by variable $s_{jk}\in\mathbb{Z}_{\geq 0}$.

Travel time from $i\in\sI$ to $j\in\sJ$ is denoted by $T_{ij}^\delta$. The travel time between demand points and charging stations incurs a cost of $C^\delta$ per time unit and is penalized in the objective function to encourage close proximity assignments. 

The expected queue for each type $k$ at $j$ is modeled with an $M$/$M$/$s_{jk}$ queue, where charging demand and supply rates are assumed to follow a Poisson process. In this model, BEBs to be recharged by type $k$ at $j$ form into a single queue, and the service is provided with a first come, first served basis. Out of the $M$/$M$/$s_{jk}$ queue performance metrics, we are particularly interested in the expected wait time $\mathbb{W}(\sum_{i\in\sI_j} \lambda_i x_{ijk}, \mu_k, s_{jk})$ inclusive of queueing and charging at charger type $k$ of station $j$, and the metric is a function of demand rate serviced by the station, $\mu_k$, and $s_{jk}$. Charger utilization rate $\rho_{jk}=\sum_{i\in\sI_j} \lambda_i x_{ijk}/ (\mu_k s_{jk})$. The condition $\rho_{jk}<1$ needs to be satisfied for queue stability, and we transform this into the linear constraint form $\mu_k s_{jk} (1-\epsilon) \geq \sum_{i\in\sI_j} \lambda_i x_{ijk}~\forall j\in\sJ, k\in\sK$. The probability that all $s_{jk}$ chargers are busy is denoted by $\mathbb{P}(\sum_{i\in\sI_j} \lambda_i x_{ijk}, \mu_k, s_{jk})$ and feeds into $\mathbb{W}(\cdot)$. Closed form formulations for both terms provided below can be found in the queueing theory literature. Equations \labelcref{Prob} and \labelcref{W} adapted into our notation conventions are from \cite{winston}. We assume $\mathbb{P}(\cdot)=0$ when $s_{jk}=0$. The expected wait time to receive charging service with charger type $k$ at station $j$ is denoted by variable $W_{jk}\in\mathbb{R}_{\geq 0}$, and its associated probabilistic queue waiting time and charging service time cost is $C^\tau$ per time unit.

\begin{equation}\label[eq]{Prob}
    \mathbb{P}\left(\sum_{i\in\sI_j} \lambda_i x_{ijk}, \mu_k, s_{jk}\right) = \frac{(\rho_{jk} s_{jk} )^{s_{jk}}}{(1-\rho_{jk}) s_{jk}!} \bigg(\frac{(\rho_{jk} s_{jk})^{s_{jk}}}{(1-\rho_{jk}) s_{jk}!} + \sum_{r=0}^{s_{jk}-1}\frac{(\rho_{jk} s_{jk})^r}{r!} \bigg)^{-1}
\end{equation}
\begin{equation}\label[eq]{W}
    \mathbb{W}\left(\sum_{i\in\sI_j} \lambda_i x_{ijk}, \mu_k, s_{jk}\right) = \frac{\mathbb{P}(\sum_{i\in\sI_j} \lambda_i x_{ijk}, \mu_k, s_{jk})}{\mu_k s_{jk}(1-\rho_{jk})} + \frac{1}{\mu_k}
\end{equation}

To this end, we model the problem with an MINLP. The formulation was introduced as a classical approach to model the SLCIS-BO in \cite{bermanandkrass}. Our contribution in this general form is the inclusion of charger type $k$ and the reduction of coverage constraints using subset $\sJ_i$, as well as removing the proximity constraints that lead to better solutions.

\begin{equation}\label[objfunc]{obj_fun}
    \min \mathbb{C} = \sum_{j\in\sJ} C^\phi_j y_j + \sum_{j\in\sJ, k\in\sK} C^\xi_k s_{jk} + \sum_{i\in\sI, j\in\sJ_i, k\in\sK}\lambda_i (C^\delta T^\delta_{ij} + C^\tau W_{jk}) x_{ijk}
\end{equation}
\noindent subject to,
\begin{equation}\label[consset]{active_if_visited}
    x_{ijk} \leq y_j \qquad \forall i\in\sI_j, j\in\sJ, k\in\sK
\end{equation}
\begin{equation}\label[consset]{demand_served_by_only_one}
    \sum_{j\in\sJ_i, k\in\sK} x_{ijk} =1 \qquad \forall i\in\sI
\end{equation}
\begin{equation}\label[consset]{rho}
    \mu_k s_{jk} (1-\epsilon) \geq \sum_{i\in\sI_j} \lambda_i x_{ijk} \qquad \forall j\in\sJ, k\in\sK
\end{equation}
\begin{equation}\label[consset]{W_def}
    W_{jk} \geq \mathbb{W}(\sum_{i\in\sI_j} \lambda_i x_{ijk}, \mu_k, s_{jk}) \qquad \forall j\in\sJ, k\in\sK
\end{equation}
\begin{equation}\label[consset]{non-neg1}
    x_{ijk},~y_j \in \{0,1\} \qquad \forall i\in\sI_j, j\in\sJ, k\in\sK
\end{equation}
\begin{equation}\label[consset]{non-neg2}
    W_{jk}\in \mathbb{R}_{\geq 0},~s_{jk}\in\mathbb{Z}_{\geq 0} \qquad \forall j\in\sJ, k\in\sK
\end{equation}

Objective function \labelcref{obj_fun} minimizes the total station, charger, travel, and waiting (inclusive of charging time) costs per time unit denoted by $\mathbb{C}$. Constraints \labelcref{active_if_visited} ensure that demand point \( i \in \sI \) can use charger type \( k \in \sK \) at station \( j \in \sJ \) only if station \( j \) is active. Constraints \labelcref{demand_served_by_only_one}, called \emph{single sourcing} constraints by \cite{bermanandkrass}, enforce that each demand point \( i \in \sI \) uses exactly one charger type \( k \in \sK \) at one station \( j \in \sJ \). Constraints \labelcref{rho} and \labelcref{W_def} are related to queueing theory and set an upper bound for utilization rate and a lower bound for service inclusive waiting time, respectively. Constraints \labelcref{rho} satisfy the queue stability condition $\rho_{jk} < 1$, that is the service rate per charger type is strictly greater than the charging demand rate. Note that constraints \labelcref{W_def} can be written in equality. Constraints \labelcref{non-neg1} and \labelcref{non-neg2} define variable domains. Constraint \labelcref{non-neg1} specifies that \( x_{ijk} \) and \( y_j \) are binary variables: \( x_{ijk} = 1 \) if demand point \( i \in \sI \) uses charger type \( k \in \sK \) at station \( j \in \sJ_i \), and \( y_j = 1 \) if candidate location \( j \in \sJ \) is an active charging station; otherwise, they are 0. Constraint \labelcref{non-neg2} states that the expected charging and waiting time for charger \( k \in \sK \) at station \( j \in \sJ \) is a non-negative continuous variable, and \( s_{jk} \), the number of chargers of type \( k \in \sK \) at each charging station \( j \in \sJ \), is a non-negative integer.

In addition to the constraints presented, the SLCIS-BO literature (e.g., \citet{aboolian2008location,bermanandkrass}) employs constraints \labelcref{assign_to_closest_open} called proximity constraints. The parameter $\bigM$ is a constant big number. They ensure that a demand point is assigned to the closest active station and are adopted to reduce the problem complexity. By assigning a demand point to the closest active station, we miss a potential reduction in the objective value that could be gained by a decrease in the waiting time cost. That is, assigning a particular demand point to a station that is not the closest could decrease the total system cost because the reduction in the waiting costs could outweigh the increase in the travel costs.

\begin{equation}\label[consset]{assign_to_closest_open}
    \sum_{j^\prime\in\sJ_i, k\in\sK} T^\delta_{ij^\prime} x_{ij^\prime k} \leq (T^\delta_{ij} - \bigM) y_{j} + \bigM \qquad \forall i\in\sI_j, j\in\sJ
\end{equation}

Compared to the model in \citet{aboolian2008location}, our model targets a larger problem because we do not constrain each $i \in \sI$ to be assigned to the closest (shortest travel time) active $j \in \sJ$. In our model, all decisions are driven by $x_{ijk}$ variables. For a given vector $\textbf{x}$, we can find the best solutions to all other variable vectors. When we have the optimal $\textbf{x}$, we can prove optimality. Yet, there are $2^{\norm{\{(i,j) | i\in\sI\land j\in\sJ_i \}} \norm{\sK}}$ possible different solutions to $\textbf{x}$. In their study, the authors assumed that, dropping index $k \in \sK$ without loss of generality, $x_{ij}=1$ for $i\in\sI$ with the closest active station $j\in\sJ$. With such an approach, the problem complexity reduces to $2^{\norm{\sJ}}$ because the best solutions to all variables can be derived once $\textbf{y}$ is given. Even if $\textbf{y}$ is given in our model, we cannot infer the best solutions to $\textbf{x}$ because the model makes $x_{ijk}$ decisions not only based on the travel time between $i\in\sI$ and $j\in\sJ$ but also considering the queueing time and the particular $k\in\sK$ used. For this reason, the exact solution in \citep{aboolian2008location} is not directly applicable to our model although our model considerably resembles the model in the literature. Regarding the difficulty introduced by removing the proximity constraints, \citet{berman2007multiple} stated ``\emph{Such a model is also much more complicated, and we can only attempt to find an equilibrium if one exists}''.

\section{Solution approaches} \label[sec]{approaches}

\subsection{An exact solution method}\label[sec]{exact_method}

\citet{grassmann1983convexity} proved that $\sum_{i\in\sI_j}\lambda_i x_{ijk}\mathbb{W}(\cdot)$ is convex in $\rho_{jk}$. \citet{lee1983note} proved the Erlang delay formula to be strictly increasing and convex in $\rho_{jk}$, and noted that $\mathbb{W}(\cdot)$ satisfies the same condition as in \cref{W_convex_in_rho}. This property is used to derive an exact solution method that relies on the cutting-plane method.

\begin{property}\label[proper]{W_convex_in_rho}
    $\mathbb{W}(\cdot)$ is strictly increasing and convex in $0<\rho_{jk}<1$.
\end{property}

The model is nonlinear due to $W_{jk}$. Let $q_{ijk}\in\mathbb{R}_{\geq 0}$ denote the travel and charging time inclusive waiting costs per time unit, which is the last component of the objective function \labelcref{obj_fun}. Let $\overline{S}_{jk}$ be the maximum number of chargers of type $k$ that can be allocated to $j$, and let $\sS=\{0,1,2,\ldots,\max_{j\in\sJ,k\in\sK}\overline{S}_{jk}\}$ denote a set of all possible discrete number of chargers to be allocated. The number of chargers given $c\in\sS$ is then denoted by the parameter $S_c$. We introduce the binary variable $z_{cjk}=1~\forall c\in\sS,~j\in\sJ,~k\in\sK$ to indicate if $c$ number of type $k$ chargers are allocated to $j$. We now need to restrain $z_{cjk}=1$ for only one $c\in\sS$ given $j$ and $k$. Finally, we replace \labelcref{obj_fun} with \labelcref{linear_obj_fun} and introduce \labelcref{linear_obj_fun2}--\labelcref{non-neg4}. After these definitions and replacements, $W_{jk}$ and therefore constraints \labelcref{W_def} and \labelcref{linear_obj_fun2} remain non-linear, which are addressed next.

\begin{equation}\label[objfunc]{linear_obj_fun}
    \min \mathbb{C} = \sum_{j\in\sJ} C^\phi_j y_j + \sum_{j\in\sJ, k\in\sK} C^\xi_k s_{jk} + \sum_{i\in\sI, j\in\sJ_i, k\in\sK}q_{ijk}
\end{equation}
\begin{equation}\label[consset]{linear_obj_fun2}
    q_{ijk} \geq \lambda_i (C^\delta T^\delta_{ij}+C^\tau W_{jk}) - \bigM(1-x_{ijk}) \qquad \forall i\in\sI, j\in\sJ_i, k\in\sK
\end{equation}
\begin{equation}\label[consset]{s_to_z}
    s_{jk} = \sum_{c\in\sS} S_c z_{cjk} \qquad \forall j\in\sJ, k\in\sK
\end{equation}
\begin{equation}\label[consset]{s_to_z2}
    \sum_{c\in\sS} z_{cjk} = 1 \qquad \forall j\in\sJ, k\in\sK
\end{equation}
\begin{equation}\label[consset]{non-neg3}
    q_{ijk}\in \mathbb{R}_{\geq 0} \qquad \forall i\in\sI_j, j\in\sJ, k\in\sK
\end{equation}
\begin{equation}\label[consset]{non-neg4}
    z_{cjk} \in \{0,1\} \qquad \forall c\in\sS, j\in\sJ, k\in\sK
\end{equation}

Following \cref{W_convex_in_rho} and the cutting-plane method proposed in \cite{cokyasar2023additive}, we will iteratively enforce lower bounds to $\mathbb{W}(\cdot)$. Let $\mathbb{W}^\nu(\rho_{jk})$ defined in \labelcref{W_nu} represent a portion of $\mathbb{W}(\cdot)$, and notice that \cref{W_convex_in_rho} applies to $\mathbb{W}^\nu(\rho_{jk})$ as well. 

\begin{equation}\label[eq]{W_nu}
    \mathbb{W}^\nu(\rho_{jk}) = \frac{\mathbb{P}(\sum_{i\in\sI_j} \lambda_i x_{ijk}, \mu_k, s_{jk})}{1-\rho_{jk}}
\end{equation}

Let $\Tilde{\rho}_{jk}$ represent an approximation to $\rho_{jk}$, $B_{jk}=\partial \mathbb{W}^\nu(\Tilde{\rho}_{jk}) / \partial \Tilde{\rho}_{jk}$, and $A_{jk}=\mathbb{W}^\nu(\Tilde{\rho}_{jk}) - B_{jk} \Tilde{\rho}_{jk}$. Then, the linear line $A_{jk}+B_{jk}\rho_{jk}$ supports $\mathbb{W}^\nu(\cdot)$ at $\rho_{jk}=\Tilde{\rho}_{jk}$, and we write $W_{jk} \geq (A_{jk}+B_{jk}\rho_{jk})/(\mu_k S_c) + 1/\mu_k$. Introducing this constraint conditional upon values of $S_c$, $A_{jk}$, $B_{jk}$, $x_{ijk}$, and $z_{cjk}$ will ensure $W_{jk}$ is provided a lower bound to be greater than or equal to the true value of $\mathbb{W}(\cdot)$. We can initially remove \labelcref{W_def} from the model and introduce these lazy constraints every time an integer solution is obtained during the branch-and-bound process with \labelcref{W_cut}.

\begin{equation}\label[consset]{W_cut}
    W_{jk} \geq \frac{A_{jk} z_{S_c jk}}{\mu_k S_c} + \frac{B_{jk} \sum_{i\in \sI_j} \lambda_i (x_{ijk} + z_{S_c jk} - 1)}{\mu_k^2 S_c^2} + \frac{z_{S_c jk}}{\mu_k} \qquad \forall j\in\sJ, k\in\sK
\end{equation}

Constraints \labelcref{W_cut} ensure that $W_{jk}$ is at least the value on the right hand side of the constraint if type $k$ charger at station $j$ is deployed $S_c>0$ number of chargers. In case $S_c=0$, the constraint is undefined. The solution procedure is as follows. Build the model removing \labelcref{W_def}. Within the mixed-integer programming (MIP) solver, using the solution vectors $\mathbf{s}$, $\mathbf{x}$, and $\mathbf{y}$ every time a feasible MIP solution is obtained, i) fetch $\mathbb{C}^\iota$, the lower bound from the solver, ii) compute $\Tilde{\rho}_{jk}= \sum_{i\in\sI_j}\lambda_i x_{ijk}/(\mu_k s_{jk})$, iii) compute $\overline{W}_{jk}$, the upper bound (i.e., true value) of $W_{jk}$ via $\Tilde{\rho}_{jk}$, iv) compute $\mathbb{\overline{C}}$, the upper bound (i.e., true value) of $\mathbb{C}$ via $\overline{W}_{jk}$, iv) the solution gap $Q:=1-(\mathbb{C}^\iota/\mathbb{\overline{C}})$. While $Q>Q^\tau$, where $Q^\tau$ is an acceptable gap threshold, let $S_c=s_{jk}$ and introduce \labelcref{W_cut} to $j$ and $k$ only when $s_{jk}>0$ and $W_{jk} < \overline{W}_{jk}$. Therefore, constraints \labelcref{W_cut} are only introduced to specific charger types at stations that the current solution shows that i) there is at least one charger allocated, ii) the value of $W_{jk}$ is underestimated. Once $Q\leq Q^\tau$ is satisfied (or a predefined solution time limit is reached), the solver terminates and reports the best solution found. Notice that in the second step of the procedure, the solver is not guaranteed to provide the true value of $W_{jk}$  unless a lower bound is provided, and the only initial bounding condition for the variable is $W_{jk}\geq 0$. Allowing the solver exhaustively introduce these lazy constraints will ensure reaching to an optimal solution. We use Gurobi's MIP callback \citep{gurobicallback} functionality in the programmatic application.

The exact solution method presented in this section can handle the model with and without constraints \labelcref{assign_to_closest_open}. However, the method, similar to many exact solution methods, does not provide accurate results as the problem size grows. To this end, we adopt metaheuristic solution algorithms similar to those in \citep{berman2007multiple}.

\subsection{Simulated Annealing}\label[sec]{SA}
We begin with providing definitions and properties that shape the algorithmic procedure in this section.

\begin{definition}\label[def]{minsjk_def}
    Let $\sX$ denote a solution set of $(i,j,k)$ triplets whose $x_{ijk}=1$, that is $\sX=\{(i,j,k)|x_{ijk}=1\land i\in\sI_j\land j\in\sJ \land k\in\sK\}$, and let $s_{jk}^{\min}$ denote the minimum number of type $k\in\sK$ chargers needed at station $j\in\sJ$ to satisfy constraints \labelcref{rho}, that is $s_{jk}^{\min}=\sum_{i^\prime\in\sI_j|i=i^\prime} \frac{\lambda_{i^\prime}}{\mu_k}~\forall (i,j,k)\in\sX$.
\end{definition}

\begin{definition}\label[def]{W_decrease_def}
    Let $\mathbb{W}(\cdot)^{+1}:=\mathbb{W}(\sum_{i\in\sI_j} \lambda_i x_{ijk}, \mu_k, s_{jk})$ – $\mathbb{W}(\sum_{i\in\sI_j} \lambda_i x_{ijk}, \mu_k, s_{jk}+1)$.
\end{definition}

\begin{property}\label[proper]{W_decreasing_in_k}
    $\mathbb{W}(\cdot)$ is decreasing in $s_{jk}$.
\end{property}

\begin{property}\label[proper]{W_convex_in_k}
    $\mathbb{W}(\cdot)$ is convex in $s_{jk}$, that is $\mathbb{W}(\cdot)^{+1}$ is decreasing in $s_{jk}$.
\end{property}

\cref{minstations} finds the minimum number of active stations ($\norm{\sJ^{\min}}$) that are necessary to have a feasible solution. In the algorithm, candidate stations in $\sJ$ are activated such that each covers the maximum number of demand points. While $\norm{\sJ^{\min}}$ is exact, the elements in $\sJ^{\min}$ may vary depending on the order of the algorithmic execution. However, the subset is useful to initialize a heuristic or metaheuristic algorithms. The solution methods in \citet{berman2007multiple}, for instance, considered a starting point of activating one $j$. That solution was feasible for their problem because any $j\in\sJ$ was a potential location for every $i\in\sI$, and impractical assignments were eliminated through coverage constraints. The subset representation of the coverage constraints enables us to use \cref{minstations} and initialize solution methods with feasible and potentially higher quality solutions.

\defrefstwo{minsjk_def}{W_decrease_def} and \proprefstwo{W_decreasing_in_k}{W_convex_in_k} lead to \cref{s_given_x} that finds the best $s_{jk}$ denoted by $s_{jk}^\prime$ for given $\sX$. In \cref{minsjk_def}, we define the minimum number of type $k$ chargers at station $j$ needed for given $\sX$ to satisfy constraints \labelcref{rho}. However, $s_{jk}^{\min}$ is not guaranteed to be the best number for the associated $\sX$ because increasing $s_{jk}^{\min}$ could drop $\mathbb{W}(\cdot)$ and hence the objective function $\mathbb{C}$. With \cref{s_given_x}, we show that the best solution to $\mathbf{s}$ vector can be obtained for $\sX$. Therefore, if $\mathbf{x}=\mathbf{x}^*$, $\mathbf{C}^*$ is the optimal solution, and associated values to variables of the optimal solution can be calculated. The algorithm begins with finding $s_{jk}^{\min}$ for tuple $(j,k)$, and then increments $s_{jk}^{\min}$ by one until no improvement in the objective function is made and incrementing is feasible. Finally, it reports the best number of type $k$ chargers to be allocated to station $j$ considering $i$ to $(j,k)$ coupling in $\sX$. This algorithm is useful because one can enumerate all applicable combinations of $\sX$ and find the optimal solution if exists.

\cref{x_calc} assigns demand points \(i \in \sI\) to an active station \(j \in \sJ^\alpha_i\) that can serve $i$. The algorithm takes as input \(\sJ^\alpha\), the set of active stations where \(y_j = 1 ~\forall j \in \sJ\), and a parameter \(P\), which ranges between 0 and 1. The output is the decision variable \(x_{ijk} \in \{0,1\} ~\forall i \in \sI, j \in \sJ, k \in \sK\), indicating which station and charger type a demand point is assigned to. The algorithm primarily assigns each demand point to its nearest active station, minimizing \(T_{ij}^\delta\), travel time between \(i\) and \(j\). However, considering the inclusion of waiting time at charging stations in our problem, the closest station may not always yield the optimal solution. To account for this, the parameter \(P\) introduces randomness, determining the probability of assigning demand point \(i\) to a station \(j \in \sJ^\alpha_i\) other than the closest one. Additionally, the algorithm randomly selects the charger type \(k \in \sK\).

SA is a metaheuristic method inspired by the cooling process of melted metals \citep{berman2007multiple}. We utilize three main parameters for SA: the initial temperature $T_0$, the maximum number of iterations $L$, and the temperature decreasing factor $C$. \cref{SA_algo} shows the steps of the SA algorithm for solving this problem.

At each iteration, the \textproc{SimulatedAnnealing($\cdot$)} function in \cref{SA_algo} generates a new subset of stations \(\sJ^\nu\) by either randomly closing one of the active stations in the current set \(\sJ^\theta\) or activating one of the closed stations. This process is repeated until a feasible subset \(\sJ^\nu\) is found, ensuring that the active stations cover all demand points. Using \algorefstwo{s_given_x}{x_calc}, the algorithm assigns demand points to stations and determines the optimal number of each charger type $k \in \sK$ at each facility $j \in \sJ^\nu$. The objective value of the \textit{new stations} is then evaluated, and if it improves upon the best-known solution (i.e., \(\mathbb{C}^\nu < \mathbb{C}^\beta\)), the best solution is updated accordingly. If not, the Metropolis acceptance criterion \(m\) is calculated, which determines the probability of accepting the new solution as the current solution. This probability is based on the current temperature \(T\), the current objective value \(\mathbb{C}^\theta\), and the best objective value \(\mathbb{C}^\beta\). The temperature \(T\) is then updated, and the process repeats. The algorithm continues iterating until the maximum number of iterations \(L\) is reached.

\subsection{Genetic Algorithm}\label[sec]{GA}
\cref{covstations} finds a population $\sJ^{\sC}$ of size $N$ consisting of the smallest subsets of charging stations $\sJ^\alpha \in \sJ^\sC$. Each subset can cover all demand points and will contain either the minimum number of active stations or the minimum plus one number of active stations. This set is determined using the function \textproc{CoverSets}$(\cdot)$, which employs a backtracking algorithm defined in the \textproc{BackTrack}$(\cdot)$ function.

The \textproc{BackTrack}($\sI^{\kappa}$, $\sJ^{\alpha}$, $\sJ^{\kappa}$, $\sJ^\sC$, $S$) function in \cref{covstations} systematically explores all possible combinations of activating charging stations, to find the smallest subsets that can cover all demand points. This function takes as input a set of remaining (uncovered) demand points, $\sI^\kappa$, a set of currently active stations, $\sJ^\alpha$, and a set of remaining (inactive) charging stations, $\sJ^\kappa$. It operates recursively, activating stations one by one and checking if the new set of activated stations can cover all remaining demand points. The function improves efficiency by pruning paths that cannot lead to a solution with a small subset of stations.

GAs are metaheuristic methods that simulate the evolutionary processes of natural selection and survival of the fittest \citep{YiHe2023}. The process begins with creating an initial population \(\sJ^\sC\), where each individual, or \textit{chromosome}, represents a feasible configuration of open charging stations \(\sJ^\alpha \in \sJ^\sC\). The fitness of each chromosome is evaluated via \( \mathbb{C} \). Through selection, crossover, and mutation operations, new chromosomes are generated in each iteration. This process continues until the termination criterion, defined by reaching the maximum number of iterations \(L\), is met, resulting in the best solution found.

In this algorithm, chromosomes represent combinations of open stations in the population \(\sJ^\alpha \in \sJ^\sC\), ensuring that all \(i \in \sI\) are covered. After generating the initial population \(\sJ^\sC\) using \cref{covstations}, the \textit{selection} phase begins. A subset \(\sJ^\sP \subseteq \sJ^\sC\) is randomly selected from the population, and two chromosomes with the highest fitness (lowest \( \mathbb{C} \)) in this subset are chosen as parents (\(\sJ^{\pi1} \And \sJ^{\pi2}\)). The parameter \(Y\) specifies the size of \(\sJ^\sP\). During the \textit{crossover} phase, the offspring chromosome inherits genes from parent 1 for the first half of the candidate stations (\( y^\nu_j \coloneqq y^{\pi 1}_j \)) and from parent 2 for the second half (\( y^\nu_j \coloneqq y^{\pi 2}_j \)). Since each charging station \(j \in \sJ\) serves a specific \(i \in \sI_j\), we ensure that the offspring’s active stations \(\sJ^\nu\) cover all demand points. If any demand point is left uncovered, random candidate locations in \(\sJ^\nu\) are activated until full coverage is achieved.

Like the SA approach, the GA uses \algorefstwo{s_given_x}{x_calc} to assign demand points to active stations \(\sX^\nu\) and to determine the optimal number of each charger type \(s_{jk}^\nu ~\forall k \in \sK\) at each facility \(j \in \sJ^\nu\). However, in GA, the offspring inherits charger types \(k \in \sK\) based on the parental configuration, ensuring consistency in charger type allocation. A mutation step introduced in \cref{x_calc}, which randomly assigns some demand points to stations other than the nearest one, allows the algorithm to explore alternative solutions without being constrained to assigning demand points to the nearest open station \labelcref{assign_to_closest_open}.

Finally, the GA evaluates the fitness of the offspring \(\mathbb{C}^\nu\) by calculating the objective function as outlined in Equation \labelcref{obj_fun}. If the offspring's objective value is better than the best value in the population (i.e., \(\mathbb{C}^\nu<\mathbb{C}^\beta\)), the best solution is updated \(\mathbb{C}^\beta \coloneqq \mathbb{C}^\nu\), and the offspring $\sJ^\nu$ replaces the worst solution $\sJ^\omega$ in the population. If the offspring is neither better than the best solution nor worse than the worst (i.e., \(\mathbb{C}^\beta \leq \mathbb{C}^\nu \leq \mathbb{C}^\omega\)), it still replaces $\sJ^\omega$. However, if it is worse than the worst solution (i.e., \(\mathbb{C}^\nu > \mathbb{C}^\omega\)), \(\sJ^\omega\) is replaced with \(\sJ^\nu\) based on a probability of \(\frac{\mathbb{C}^\omega}{\mathbb{C}^\nu}\). \cref{fig_GA} presents the flowchart of the GA implemented in this study.

\section{Numerical experiments} \label[sec]{num_exp}
In this section, we demonstrate the findings of extensive numerical experiments. We first elaborate the design of experiments by referencing to data sources adopted. The results of computational experiments that compare the three solution methods from solution quality and solution time perspectives are presented. Case studies focusing on CTA and Pace agencies illustrate the importance of collaboration in the electrification context, as well as the gain of considering locations other than garages as candidate charging stations. Finally, sensitivity analyses are conducted to identify the impact of key modeling parameters on the overall system cost. 

\subsection{Design of experiments} \label[sec]{exp_design}
Various data sources are integrated to design the numerical experiments for deploying charging stations in the Chicago metropolitan region in Illinois, USA. The data related to trip schedules are publicly available through the General Transit Feed Specification (GTFS) \citep{GTFS}. We obtained the data for 25,999 trips that are operated by CTA and Pace. These trips use 1,078 bus stops as origins and destinations. Using the model in \cite{ADavatgari2024} and implementing the logic described in the first step, we identified 414 stops where the charging need arises. With regard to the candidate charging stations, we obtained 30,236 stops used by two agencies, seven CTA bus garages, and ten Pace bus garages. Therefore, 414 stop locations and 30,253 stop locations, depicted in \cref{study_region}, are considered as a pool of demand points and candidate charging stations, respectively. These locations are then partitioned into clusters to control the problem size and to conduct the analyses in this section. Two charger types, fast and slow, are considered to be allocated. Parametric details related to schedule formation are provided in \cite{ADavatgari2024}. \cref{param_sum} summarizes the parameters used in the baseline scenario.

According to \cite{CTAAgencyProfile}, CTA buses accumulate 4,830,866 annual revenue hours, with operating expenses totaling \$774,665,363. Hence, the hourly cost of a travel is determined by dividing the annual operating expenses by the annual revenue hours. For charging and waiting time, an additional cost related to energy, obtained from \cite{AvgEnergyPrice}, is added. Specifically, a charge of \$0.79 per minute supplements the travel cost.

We assume $Q^\iota=10$ (\%) and $Q^\phi=80$ (\%). Time to recharge BEB \( T_k^\rho \) is then calculated as follows: 
$T_k^\rho = \frac{B \times 3,600 \times (Q^\phi- Q^\iota)}{P_k}$
where \( B \) is in kWh, and \( P_k \) is in kW. The service rate \( \mu_k \) is then $\mu_k= (T_k^\rho)^{-1}$. All experiments are conducted converting time units to minutes and cost units to USD per minute.

\subsection{Computational experiments} \label[sec]{comp_exp}
We compare the computational performance of the exact method (henceforth called MINLP) solved through Gurobi \citep{gurobi}, SA, and GA. This comparison allows us to assess their scalability and to identify the solution quality they provide. We begin with clustering 414 stops into 100 charging demand points and clustering 30,253 stops into 100 candidate charger locations, that is we have $\norm{\sI}=100$ and $\norm{\sJ}=100$. We then sample scenarios with different sizes of $\sI$ and $\sJ$ using this pool. We consider two main scenarios: i) $\norm{\sJ}=5$ and $\norm{\sI}\in [3,5,7,10,20,50]$, ii) $\norm{\sI}=5$ and $\norm{\sJ}\in [3,5,7,10,20,50]$. The goal with this categorization is to observe the computational quality with varying dimensions of the problem in both $\sI$ and $\sJ$. The number of instances for the smallest size problems are set to 500, and the number is decreased as problem sizes grow. Each instance for each solution method is provided with a fixed computational time limit. The limit in the smallest instances is 30 seconds and is increased in conjunction with problem sizes as seen fit. The acceptable gap threshold for the MINLP ($Q^\tau$) is set to 0.01\% for all instances. Each instance is solved ten times in parallel using both SA and GA, and the best solution is reported. In this section, a total of 1,470 instances were solved allocating 55 hours of computational time. All computations were carried out on an Intel{\textsuperscript \textregistered} Xeon{\textsuperscript \textregistered} Gold 6138 CPU @2.0 GHz workstation with 128 GB of RAM and 40 cores. Problem instances were solved by using the Python 3.12.6 interface to the commercial solver Gurobi 9.1.2~\citep{gurobi}. 

\begin{table}[!htbp]
\centering
\caption{Comparing the solution quality statistics for MINLP, SA, and GA methods.}
\label[tab]{quality_comparison}
\begin{center}
\resizebox{\textwidth}{!}{%
\begin{tabular}{ccccccccc}
\toprule
\multicolumn{4}{c}{Scenario} & \multicolumn{1}{c}{MINLP} & \multicolumn{2}{c}{SA} & \multicolumn{2}{c}{GA} \\
\cmidrule(lr){1-4} \cmidrule(lr){5-5} \cmidrule(lr){6-7} \cmidrule(lr){8-9}

$\norm{\sI}$ & $\norm{\sJ}$ & \# instances & time (s) & \# optimal & \# optimal & \# better & \# optimal & \# better \\
\midrule
3 & 5 & 500 & 30 & 469 & 469 & 18/31 & 469 & 18/31 \\
5 & 5 & 200 & 60 & 186 & 186 & 11/14 & 186 & 11/14 \\
7 & 5 & 50 & 60 & 18 & 18 & 11/32 & 18 & 11/32 \\
10 & 5 & 50 & 60 & 0 & - & 24/50 & - & 24/50 \\
20 & 5 & 25 & 120 & 0 & - & 3/25 & - & 21/25 \\
50 & 5 & 10 & 300 & 0 & - & 1/10 & - & 9/10 \\
\midrule
5 & 3 & 500 & 30 & 460 & 460 & 28/40 & 437 & 27/40 \\
5 & 7 & 50 & 60 & 35 & 35 & 8/15 & 35 & 8/15 \\
5 & 10 & 50 & 60 & 2 & 2 & 19/48 & 2 & 19/48 \\
5 & 20 & 25 & 120 & 0 & - & 14/25 & - & 14/25 \\
5 & 50 & 10 & 300 & 0 & - & 10/10 & - & 10/10 \\
\bottomrule
\end{tabular}
}
\end{center}
\end{table}

In \cref{quality_comparison}, the scenario columns tabulate the details of the problem instances. Then, we report the number of optimal solutions for each solution method. For the solutions that are not optimal through the MINLP, we report the number where SA (or GA) provided better solutions than the MINLP. For instance, in $\norm{\sI}=3$ and $\norm{\sJ}=5$ scenario, the MINLP does not find an optimal solution to 31 instances within the provided time limit of 30 seconds. Out of these 31 instances, SA and GA provided better solutions for 18 instances compared to the MINLP. When the MINLP does not find any optimal solutions for a particular scenario, \# better column is especially useful to observe the performance of the two methods compared to the MINLP. Overall, we have two main interpretations based on this table: i) SA and GA perform quite similar especially in small instances, and GA performs relatively better than SA in large instances. ii) Problems become more difficult to solve as $\norm{\sI}$ increases compared to a relative increase in $\norm{\sJ}$.

\cref{gap_comparison} complements \cref{quality_comparison} by providing statistics of optimality gap percentages for the scenarios. The gap values reported for SA and GA are calculated by fetching the lower bound reported by the MINLP solver ($\mathbb{C}^\iota$), obtaining the upper bound found from the corresponding solution method ($\mathbb{C}^\beta$), and computing $100(1-\mathbb{C}^\iota/\mathbb{C}^\beta)$. We deduce that GA outperforms the MINLP for all scenarios on average. SA performs similarly to GA in the small instances, however its quality reduces in the large instances. Overall, we observe small gap averages in the small instances, however the gaps are considerably high in large instances. Yet, this is not a sign of poor quality because the solver spends the majority of the effort in the branch-and-bound to improve the lower bound. The solver finds a relatively good upper bound in the first few minutes, and the lower bound approaches to the upper bound while the upper bound does not considerably decrease. Eventually, the solver proves an optimality certificate by spending a considerable amount of time (days in large instances) to find a feasible lower bound that is equivalent to the upper.

\begin{table}[!ht]
\centering
\caption{Comparing the solution gap for MINLP, SA, and GA methods.}
\label[tab]{gap_comparison}
\begin{center}
\resizebox{\textwidth}{!}{%
\begin{tabular}{cccccccccccccc}
\toprule
\multicolumn{2}{c}{Scenario} & \multicolumn{4}{c}{MINLP gap (\%)} & \multicolumn{4}{c}{SA gap (\%)} & \multicolumn{4}{c}{GA gap (\%)} \\
\cmidrule(lr){1-2} \cmidrule(lr){3-6} \cmidrule(lr){7-10} \cmidrule(lr){11-14}

$\norm{\sI}$ & $\norm{\sJ}$ & min. & max. & std. & avg. & min. & max. & std. & avg. & min. & max. & std. & avg. \\
\midrule
3 & 5 & 0 & 48.4 & 4.4 & 0.7 
            & 0 & 48.4 & 2.9 & 0.3 
            & 0 & 48.4 & 2.9 & 0.3 \\
5 & 5 & 0 & 29.9 & 5.1 & 1.1 
            & 0 & 26.3 & 2.9 & 0.6 
            & 0 & 26.3 & 2.9 & 0.6 \\
7 & 5 & 0 & 60.6 & 20.3 & 22.4 
            & 0 & 60.6 & 20 & 22.2 
            & 0 & 60.6 & 20 & 22.2 \\
10 & 5 & 0.8 & 81.7 & 17.9 & 55.8 
            & 0.8 & 80.3 & 17.9 & 55.3 
            & 0.8 & 80.3 & 17.9 & 55.3 \\
20 & 5 & 27.4 & 90.6 & 18.9 & 74.1 
            & 29.2 & 91.1 & 18.4 & 74.9 
            & 27 & 90.6 & 19.2 & 73.5 \\
50 & 5 & 58.3 & 97.6 & 13.8 & 89.5 
            & 66.7 & 98 & 10.9 & 91.6 
            & 52.3 & 97.2 & 15.1 & 88.6 \\
\midrule
5 & 3 & 0 & 31.7 & 4.1 & 0.8 
            & 0 & 23.5 & 2.3 & 0.4 
            & 0 & 23.5 & 2.3 & 0.4 \\
5 & 7 & 0 & 56.1 & 14.8 & 7.8 
            & 0 & 50 & 14 & 7.5 
            & 0 & 50 & 13.9 & 7.5 \\
5 & 10 & 0 & 86.9 & 19.6 & 43.4 
            & 0 & 74.5 & 18 & 41.5 
            & 0 & 74.5 & 18 & 41.5 \\
5 & 20 & 20.8 & 98.4 & 17 & 70.7 
            & 29 & 84.7 & 14.6 & 67.1 
            & 29.3 & 84.8 & 14.6 & 67.1 \\
5 & 50 & 92.1 & 97.4 & 1.8 & 95.3 
            & 76.8 & 96.7 & 5.6 & 89.5 
            & 77.9 & 96.4 & 5.3 & 89.6 \\
\bottomrule
\end{tabular}
}
\end{center}
\footnotesize
\vspace{-10pt}
\emph{Note}: min: minimum, max: maximum, std: standard deviation, avg: average
\end{table}

\begin{figure}[!ht]
    \centering
    \begin{subfigure}[t]{0.48\textwidth}
        \centering
        \includegraphics[width=\linewidth]{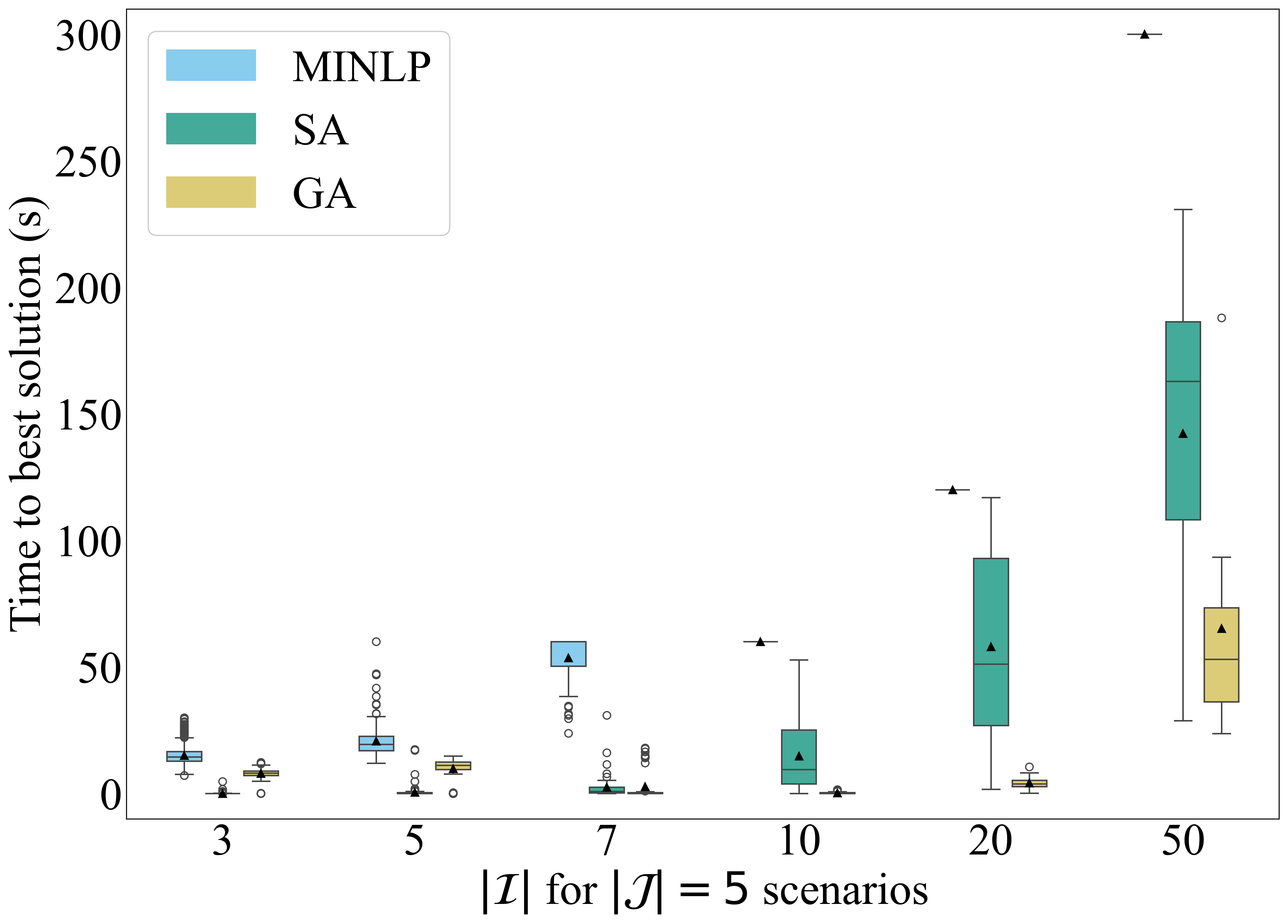} 
        \caption{Scenarios with $\norm{\sJ}=5$}
        \label[fig]{time_to_best_for_J_5}
    \end{subfigure}
    \hfill
    \begin{subfigure}[t]{0.48\textwidth}
        \centering
        \includegraphics[width=\linewidth]{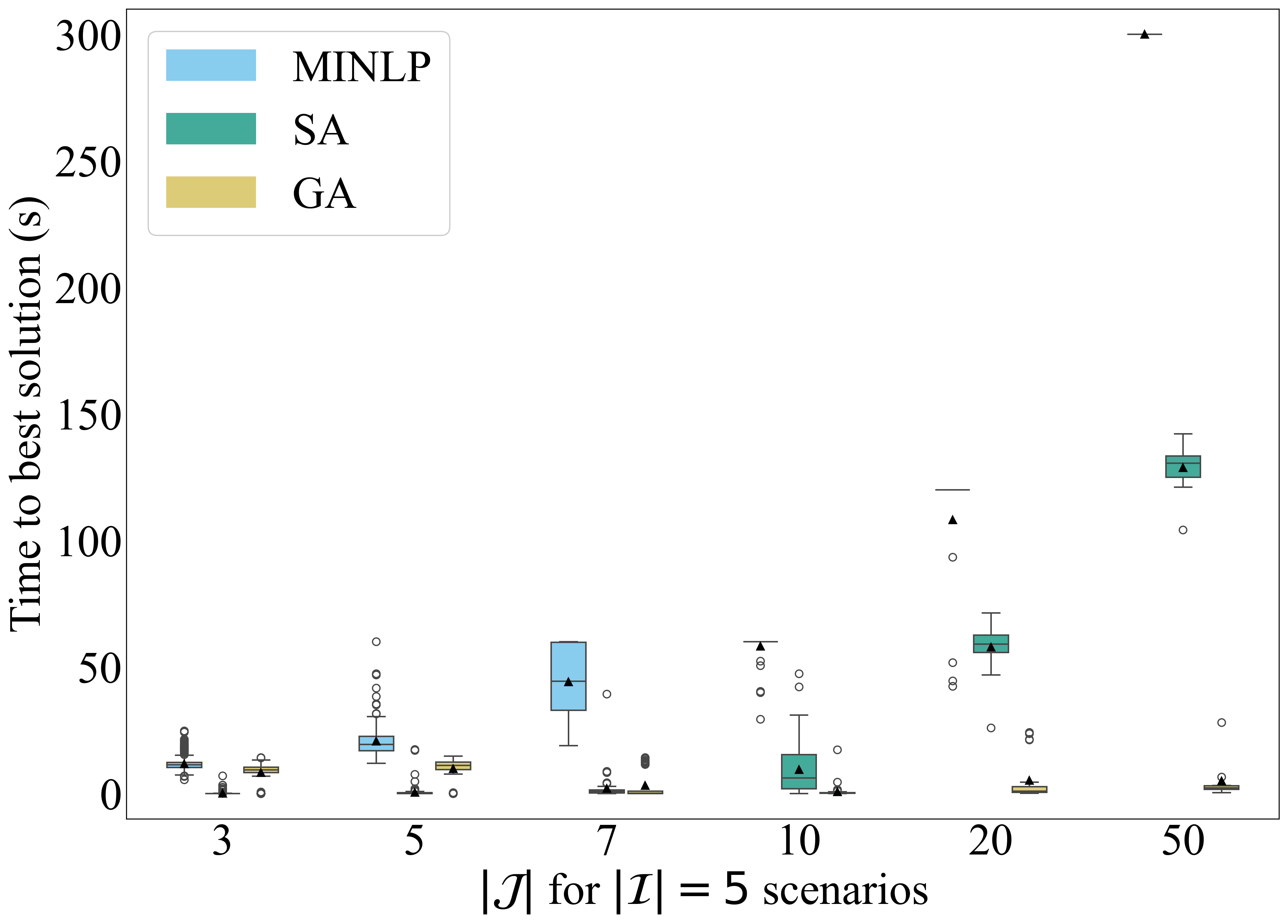} 
        \caption{Scenarios with $\norm{\sI}=5$}
        \label[fig]{time_to_best_for_I_5}
    \end{subfigure}

    \caption{Comparing the time to best solution for MINLP, SA, and GA methods.}
    \label[fig]{time_comp}
\end{figure}

The solution quality is an important factor in judging the quality of a solution method, and the time spent for finding the solution is also crucial. The problem at hand applied to practice requires patience for quality results due to the large magnitude of the real-world applications. However, one would not always possess the time, and SA and GA play a critical role in scenario analysis. \cref{time_comp} shows the time spent on finding the best solution by these three methods. SA and GA were ran the entire time during the time limit provided. Yet, we tracked the time when they found the best solution, and their statistics are reported in this figure. We can observe that SA is faster than the other two in small instances. Yet, it requires more time as the problem size increases, and we also know from the previous tables that it cannot always compete with GA. The magnitude of computational time spent on instances with $\norm{\sJ}=5$ depicted in \cref{time_to_best_for_J_5} compared to instances with $\norm{\sI}=5$ shown in \cref{time_to_best_for_I_5} also support that problem difficulty is impacted more by an increase in $\norm{\sI}$.

Since SA and GA were used ten times for each instance, the consistency in the solutions they provide is another metric to consider. \cref{SA_GA_consistency} shows the number of different solutions obtained with ten runs of SA and GA. For instance, SA provided only one single solution across ten runs in 496 instances and two different solutions across ten runs in 4 instances in the scenario with $\norm{\sI}=3$ and $\norm{\sJ}=5$. Similarly, GA had only one single solution in 499 instances and two different solutions in one instance. In the small instances, both SA and GA find consistent results. However, the consistency distorts as the scale increases. The indication is that a higher number of parallel SA and/or GA runs in large instances can yield better quality solutions.

\begin{table}[!ht]
\centering
\caption{Comparing the solution consistency statistics for SA and GA methods.}
\label[tab]{SA_GA_consistency}
\begin{center}
\begin{tabular*}\textwidth{c@{\extracolsep{\fill}}ccccccccccc}
\toprule
\multicolumn{2}{c}{Scenario} & \multicolumn{10}{c}{SA/GA}\\
\cmidrule(lr){1-2} \cmidrule(lr){3-12}
$\norm{\sI}$ & $\norm{\sJ}$ & 1 & 2 & 3 & 4 & 5 & 6 & 7 & 8 & 9 & 10\\
\midrule
3 & 5 & 496/499 & 4/1 & -/- & -/- & -/- & -/- & -/- & -/- & -/- & -/- \\
5 & 5 & 195/187 & 4/13 & 1/- & -/- & -/- & -/- & -/- & -/- & -/- & -/- \\
7 & 5 & 45/38 & 5/12 & -/- & -/- & -/- & -/- & -/- & -/- & -/- & -/- \\
10 & 5 & -/3 & 11/29 & 22/11 & 11/7 & 5/- & -/- & 1/- & -/- & -/- & -/- \\
20 & 5 & -/- & -/4 & -/5 & -/7 & -/8 & -/1 & -/- & -/- & -/- & 25/- \\
50 & 5 & -/- & -/2 & -/1 & -/3 & -/2 & -/1 & -/- & -/1 & 1/- & 9/- \\
\midrule
5 & 3 & 486/464 & 14/36 & -/- & -/- & -/- & -/- & -/- & -/- & -/- & -/- \\
5 & 7 & 48/28 & 2/20 & -/2 & -/- & -/- & -/- & -/- & -/- & -/- & -/- \\
5 & 10 & 39/5 & 9/26 & 2/13 & -/4 & -/1 & -/- & -/1 & -/- & -/- & -/- \\
5 & 20 & 11/- & 7/- & 1/1 & 2/5 & -/2 & -/7 & -/3 & -/2 & -/1 & -/- \\
5 & 50 & 4/- & 2/- & 3/- & 1/- & -/- & -/- & -/- & -/2 & -/3 & -/5 \\
\bottomrule
\end{tabular*}
\end{center}
\end{table}

Based on these tables and figures, we note that SA can be a suitable solution to small problem instances considering its time efficiency and comparable solution quality to those of MINLP and GA, and GA is a more appropriate method than SA in large problems due to its solution speed and quality.

\subsection{Case studies} \label[sec]{case_study}
The prevailing strategic approach to the deployment of BEBs in practice is housing charging equipment at bus garages. Yet, it is intuitive that, in addition to garages, deploying these devices to other locations, such as trip start and/or end terminal points, could lead to potential long-term cost savings at the cost of introducing managerial complexity, which is beyond the scope of this study. Another key challenge in the decision-making process of electrification relates to the ownership of charger station and equipment. While it is difficult to envision a scenario where all commercial sector EV users rely solely on publicly available chargers managed by third-party initiatives or local/federal authorities, this scenario is plausible for transit agencies, most of which are operated by public authorities. This is because the level of conflicting interests in public transport may be assumed lower than in the commercial sector. Consequently, charging stations for BEBs could be a joint investment among agencies operating in the same region. Such collaboration would allow multiple agencies to share charging facilities and reduce overall investment costs. 

We examine six scenarios, considering joint or separate investments, as well as charger deployment options limited to garages, other locations, or both. We use the parametric design laid out in \cref{exp_design} involving two public transit agencies in the Chicago metropolitan region as a baseline. For each scenario, we generate 16 problem instances by using problem size controllers $\sI$ and $\sJ$ and the parameters $B$, $C_k^\xi$, and $C^\tau$ as levers. In the first half of 16 instances, we group demand points and candidate charging stations into 40 clusters, that is $\norm{\sI} = 40$ and $\norm{\sJ} = 40$. In the other half, we have 50 clusters. For each half, we generate eight instances by assuming the baseline $B$ and a 50\% increased one, baseline $C_k^\xi$ and a 50\% increased one for all charger types, and a baseline $C^\tau$ and a 500\% increased one, that is two choices for each of three levers, i.e., $2^3$ instances. Eventually, we obtain 96 instances and solve each of them ten times (in parallel) using the GA method with 30 minutes of time limit for each. The results are grouped into six scenarios by taking their mean.

\cref{GAScenarios_sum} compares the results of these scenarios by showing averaged statistics, that is results in each row is the mean value of 16 problem instances. The percentage increase in the objective value \(\mathbb{C}\) relative to the baseline scenario, where agents can utilize chargers of each other, and all available candidate stations are considered, is shown in the first column following the scenario columns. The following columns respectively indicate: the average number of charging station clusters selected, the average number of slow chargers allocated, the average number of fast chargers allocated, the average of the expected waiting time across all stations in minutes, and the average utilization rates across all stations and chargers.

\begin{table}[!ht]
  \footnotesize
  \caption{Summary of scenario results.}\label[tab]{GAScenarios_sum}
  \begin{center}
   \resizebox{\textwidth}{!}{%
\begin{tabular}{ccccccccc}
\toprule
\multicolumn{3}{c}{Scenario} & \multirow{2.4}{*}{\begin{tabular}[c]{@{}c@{}}\% $\uparrow$ in $\mathbb{C}$\\ vs. baseline \end{tabular}}  & \multirow{2.4}{*}{$\sum_{j\in\sJ} y_{j}$} & \multirow{2.4}{*}{$\sum_{j\in\sJ} s_{j, \text{slow}}$} & \multirow{2.4}{*}{$\sum_{j\in\sJ} s_{j,\text{fast}}$} & \multirow{2.4}{*}{$\overline{W}_{jk} (min)$} & \multirow{2.4}{*}{$\overline{\rho}_{jk} (\%)$}\\
\cmidrule(lr){1-3}
Joint & Garage & Other & & & & & &\\ \midrule
\xmark  & \cmark  & \xmark  & 5.36  & 12.5  & 4.94  & 55.81  & 0.22  & 0.25  \\
\xmark  & \xmark  & \cmark  & 10.09 & 16.31  & 11.31  & 62.56  & 0.25  & 0.22  \\
\xmark  & \cmark  & \cmark  & 2.2 & 23.25  & 2.56  & 73.31  & 0.3  & 0.2  \\
\cmark  & \cmark  & \xmark  & 1.96 & 12  & 0.25  & 57.75  & 0.2  & 0.25  \\
\cmark  & \xmark  & \cmark  & 7.05 & 14.44  & 7.06  & 59.88  & 0.26  & 0.24  \\
\cmark  & \cmark  & \cmark  & - & 23  & 0.38  & 73.31  & 0.3  & 0.2  \\
\bottomrule
\end{tabular}}
\end{center}
\end{table}
The results show that a joint charging infrastructure deployment strategy compared to a separate one reduces the cost considerably. While garage-only scenarios lead to an increase in the cost of deployment, other-only scenarios impact the cost even worse. This highlights that transit agencies, in practice, target a more economical deployment strategy by prioritizing garage deployment and delaying charging infrastructure deployment to elsewhere. On the other hand, the solutions also depict the importance of a garage and other mixed deployment strategy. Briefly, joint and mixed deployment is the key to achieve a cost-friendly deployment of BEB infrastructure.

\subsection{Sensitivity analyses}
We performed sensitivity analyses by varying key model parameters to observe their impact in the system cost ($\mathbb{C}$). The key parameters considered are waiting time cost ($C^\tau$), charger powers ($P_k$), station cost ($C^\phi_j$), and charger cost ($C^\xi_k$). A one-at-a-time approach was used, where each parameter was individually varied while holding other parameters constant as defined in \cref{exp_design}. The results were then analyzed by comparing the objective values. 

The parameter $C^\tau$ can significantly vary from one region to another, and it could be much higher. To this end, we considered five different cases in addition to our baseline by multiplying our initial parameter by 2, 4, 6, 8, and 10. As the charging technology advances, we may expect an increase in $P_k$. Therefore, we considered four additional cases via increasing $P_k$ by 20, 40, 60, and 80 percent. Accounting for potential increases in the real estate industry, we considered five additional cases via increasing $C^\phi_j$ by 10, 30, 50, 100, and 200 percent. Finally, accounting for technological advancements that may reduce charger costs, we considered four cases via reducing $C^\xi_k$ by 20, 40, 60, and 80 percent. It is worthwhile to note that for the $P_k$ and $C^\xi_k$ parameters, each time a value was changed, it was adjusted for both types of chargers in the problem (i.e., slow and fast chargers).

A total of 36 problem instances were generated using combinations of 3, 5, 7, 10, 20, and 50 demand points and candidate station locations. For each instance, demand points and candidate stations were randomly selected from the full set of available locations (rather than clusters). 
First, we solved each instance using the baseline parametric settings earlier defined. Then, 18 other instances were created for each of 36 instances by altering one parameter at a time. In total, we solved 684 instances using GA with varying time limits based on the size of the problem. Based on $M=\norm{\sI}+\norm{\sJ}$, the time limit was set to 30, 60, 120, 300, and 420 seconds for instances with $M <10$, $10\leq M <20$, $20\leq M <50$, $50\leq M <60$, and $M\geq 60$, respectively. 

\cref{sensitivity_fig} illustrates the percentage changes in $\mathbb{C}$ in response to variations in $C^\tau$, $P_k$, $C^\phi_j$, and $C^\xi_k$. In these figures, every point corresponds to the average cost deviations of 36 instances compared to their baseline counterparts. As shown in \cref{WCost_fig}, increasing $C^\tau$ leads to a proportional increase in $\mathbb{C}$, underscoring the critical importance of minimizing waiting times during BEB recharging. Additionally, a comparison of the effects of $P_k$ and $C^\xi_k$ in \cref{ChPower_fig} and \cref{ChCost_fig} suggests that improvements in charging technology can have a greater impact if they lead to increased charger power rather than reduced charger costs. Moreover, as seen in \cref{SCost_fig}, increasing $C^\phi_j$ has a relatively minor impact on $\mathbb{C}$, with a 200\% increase in station costs resulting in only a 2\% increase in total system costs.

\section{Conclusion}\label[sec]{conclusion}
In this study, we have developed an MINLP model to tackle the PLACE-BEB effectively. Our approach began by leveraging insights from earlier studies that focused on optimizing schedules for diesel-powered buses. In the first step, we segmented these existing schedules to pinpoint potential locations where charging demand is likely to emerge, as well as suitable stations for accommodating these charging events. This segmentation process also provided the frequency of charging occurrences, which forms the charging demand.

Assuming that the charging events adhere to a Poisson process, in the second step, we adapted an SLCIS-BO model aimed at optimizing not only the locations and numbers of BEB charging stations but also the types of chargers required. During this adaptation, we identified a significant opportunity for advancement in our modeling approach by relaxing the proximity constraints. While this modification promised to enhance solution quality, it also introduced a challenge: the computational complexity escalated considerably as a result of these constraints being removed. To navigate this increased complexity and devise a viable solution methodology, we developed an exact solution approach utilizing cutting-plane techniques. However, we encountered limitations with this method, as it struggled to yield solutions for large-scale problem instances within a reasonable timeframe. Consequently, we turned our attention to the development of two metaheuristic algorithms: SA and GA.

Extensive experimental evaluations were conducted to assess the computational efficiency of the three solution methodologies. These experiments not only aimed to derive managerial insights into optimal deployment strategies through detailed case studies but also sought to understand the influence of key parameters on the overall system costs. Our findings indicated that while the exact solution method is capable of delivering optimal results for smaller instances, a one-hour computational limit was insufficient for achieving near-optimal solutions for larger instances. In contrast, SA demonstrated effectiveness for smaller cases, whereas GA outperformed the other methods when applied to large-scale instances.

Through our case studies, we highlighted the advantages of collaboratively deploying BEB chargers among various agencies operating within the same geographical region. Additionally, we considered both garage facilities and alternative locations as potential candidates for charging stations. Sensitivity analyses underscored the considerable impact of waiting time costs on minimizing overall system expenses. Therefore, we assert that the waiting time component—despite its ties to operational issues—should be integrated into the strategic planning of charger locations and allocations for electric vehicles.

While this study contributes meaningful advancements to the existing body of literature, several limitations warrant further exploration in future research. A key limitation lies in our method for determining demand points, which involved segmenting existing diesel bus schedule. This approach presents potential inaccuracies, as it fails to account for the waiting and charging times of buses during the initial demand point estimation. A discrete event simulation could validate the feasibility of the charging network provided. Based on the number of violations observed in the simulation, the network could be improved through a feedback loop between the simulation and the optimization model.

\section*{Acknowledgements}
This material is based on work supported by the U.S. Department of Energy, Office of Science, under contract number DE-AC02-06CH11357. This report and the work described were sponsored by the U.S. Department of Energy (DOE) Vehicle Technologies Office (VTO) under the Transportation Systems and Mobility Tools Core Maintenance/Pathways to Net-Zero Regional Mobility, an initiative of the Energy Efficient Mobility Systems (EEMS) Program. Erin Boyd, a DOE Office of Energy Efficiency and Renewable Energy (EERE) manager, played an important role in establishing the project concept, advancing implementation, and providing guidance.

\section*{CRediT authorship contribution statement}
\textbf{Sadjad Bazarnovi:} Formal analysis, Methodology, Software, Visualization, Writing - original draft. \textbf{Taner Cokyasar:} Conceptualization, Data curation, Formal analysis, Funding acquisition, Methodology, Software, Supervision, Visualization, Writing - original draft. \textbf{Omer Verbas:} Supervision. \textbf{Abolfazl (Kouros) Mohammadian:} Supervision.




\clearpage
\bibliographystyle{elsarticle-harv} 
\bibliography{our_bib}

\vfill
\framebox{\parbox{.90\linewidth}{\scriptsize The submitted manuscript has been created by
        UChicago Argonne, LLC, Operator of Argonne National Laboratory (``Argonne'').
        Argonne, a U.S.\ Department of Energy Office of Science laboratory, is operated
        under Contract No.\ DE-AC02-06CH11357.  The U.S.\ Government retains for itself,
        and others acting on its behalf, a paid-up nonexclusive, irrevocable worldwide
        license in said article to reproduce, prepare derivative works, distribute
        copies to the public, and perform publicly and display publicly, by or on
        behalf of the Government.  The Department of Energy will provide public access
        to these results of federally sponsored research in accordance with the DOE
        Public Access Plan \url{http://energy.gov/downloads/doe-public-access-plan}.}}

\newpage 
\renewcommand{\thesection}{Supplementary Materials}

\renewcommand{\thetable}{S.M.\arabic{table}}
\renewcommand{\thefigure}{S.M.\arabic{figure}}
\renewcommand{\thealgocf}{S.M.\arabic{algocf}}
\setcounter{table}{0}
\setcounter{figure}{0}
\setcounter{algocf}{0}

\section{Figures}\label[sup]{figures}

\begin{figure}[!ht]
    \centering
    \includegraphics[width=\linewidth]{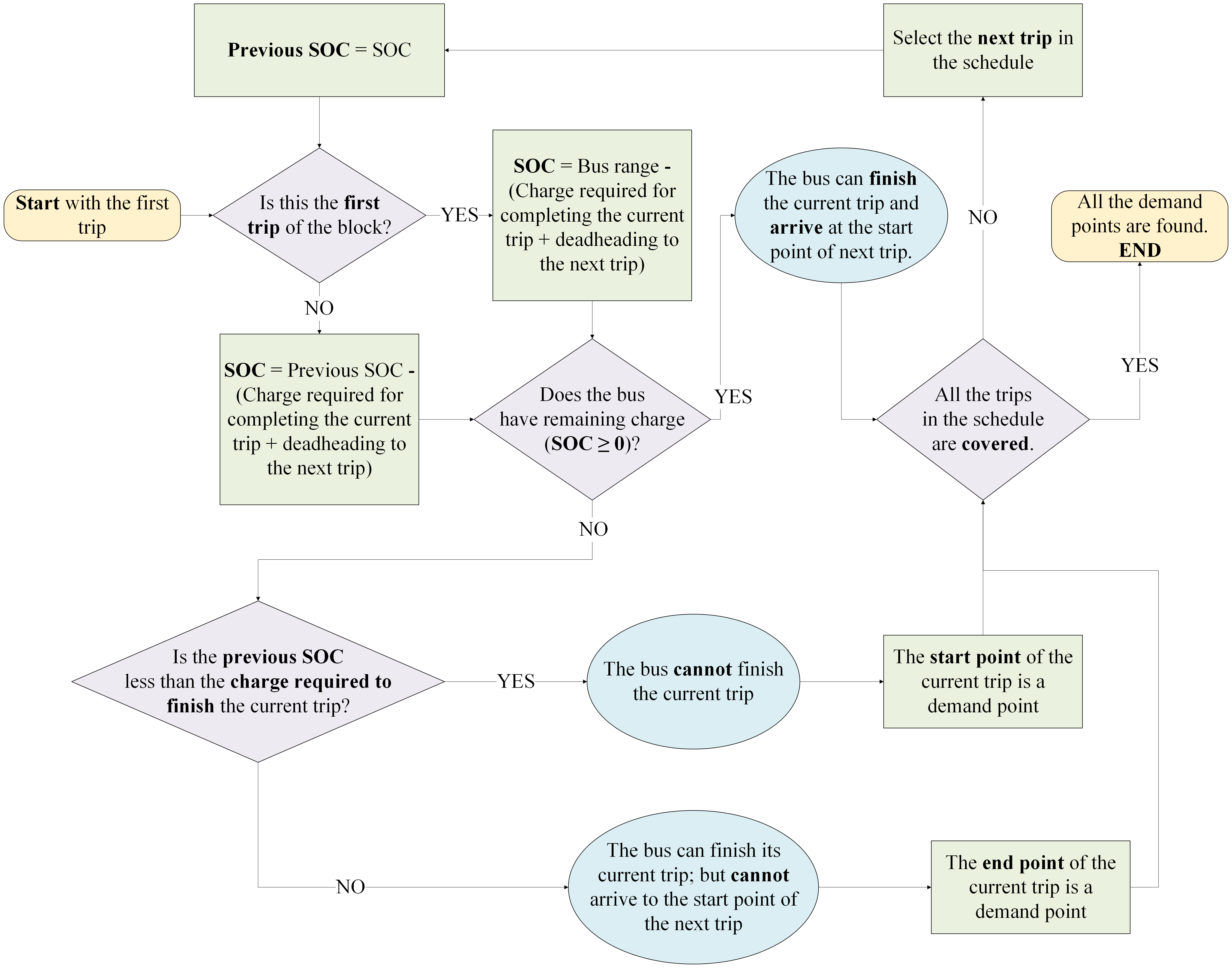}
    \caption{Flowchart for finding the demand points.}
    \label[fig]{DemandFlowChart}
\end{figure}

\begin{figure}[!ht]
    \centering
    \includegraphics[width=1\linewidth]{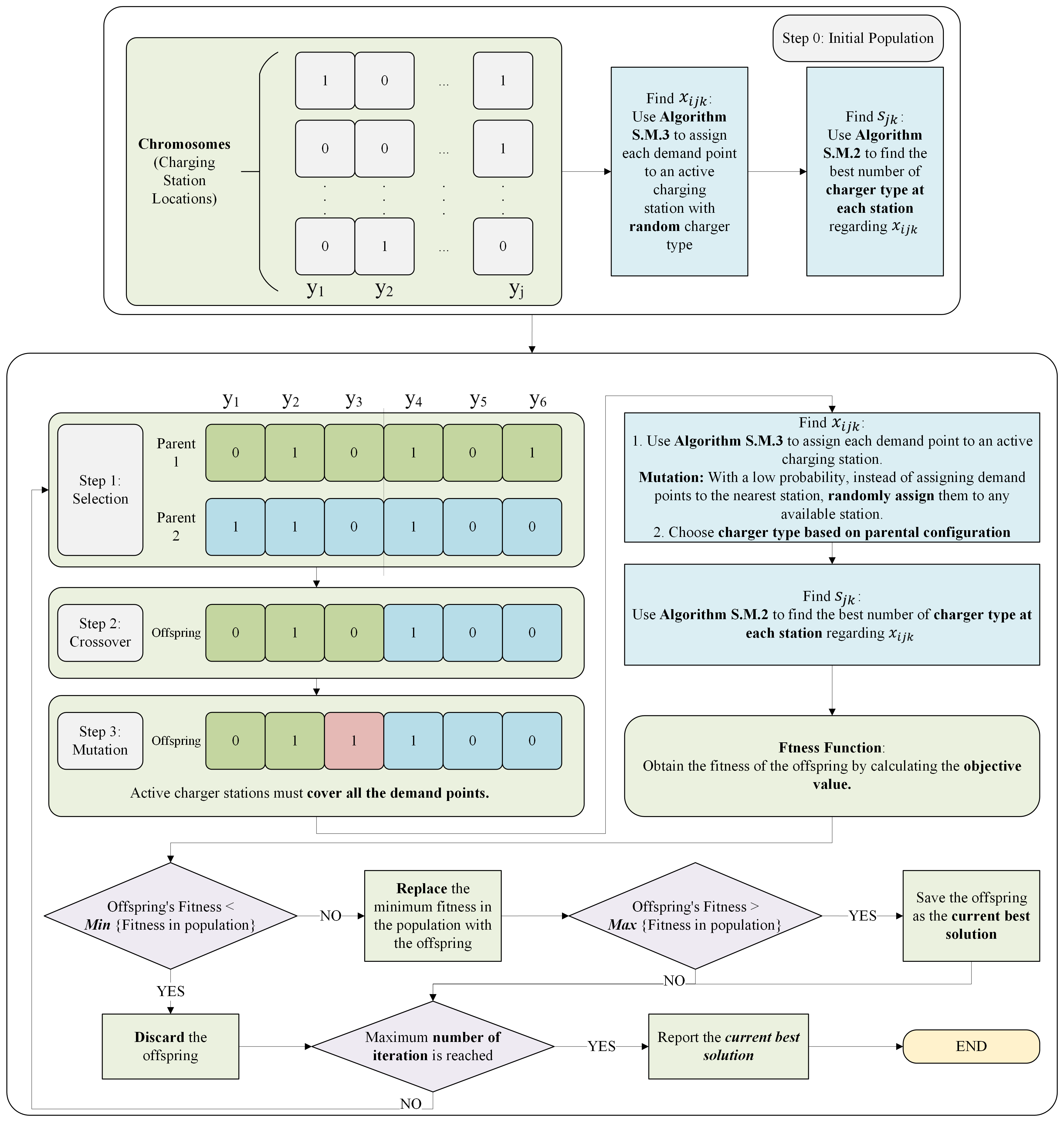}
    \caption{The genetic algorithm layout.}
    \label[fig]{fig_GA}
\end{figure}

\begin{figure}[!ht]
    \centering
    \begin{subfigure}[t]{0.49\textwidth}
        \centering
        \includegraphics[width=\linewidth]{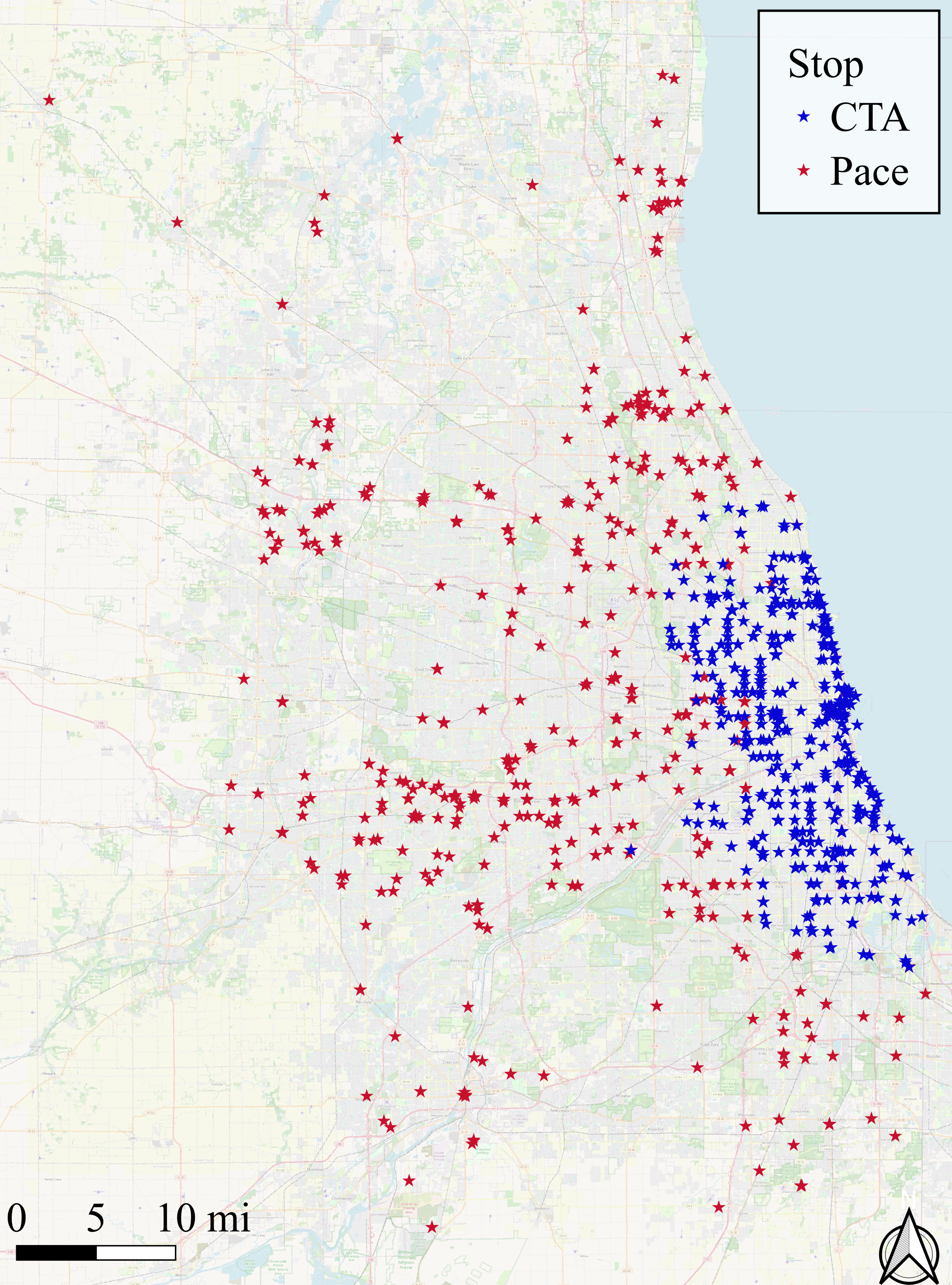} 
        \caption{Demand points}
        \label[fig]{demand_points}
    \end{subfigure}
    \hfill
    \begin{subfigure}[t]{0.49\textwidth}
        \centering
        \includegraphics[width=\linewidth]{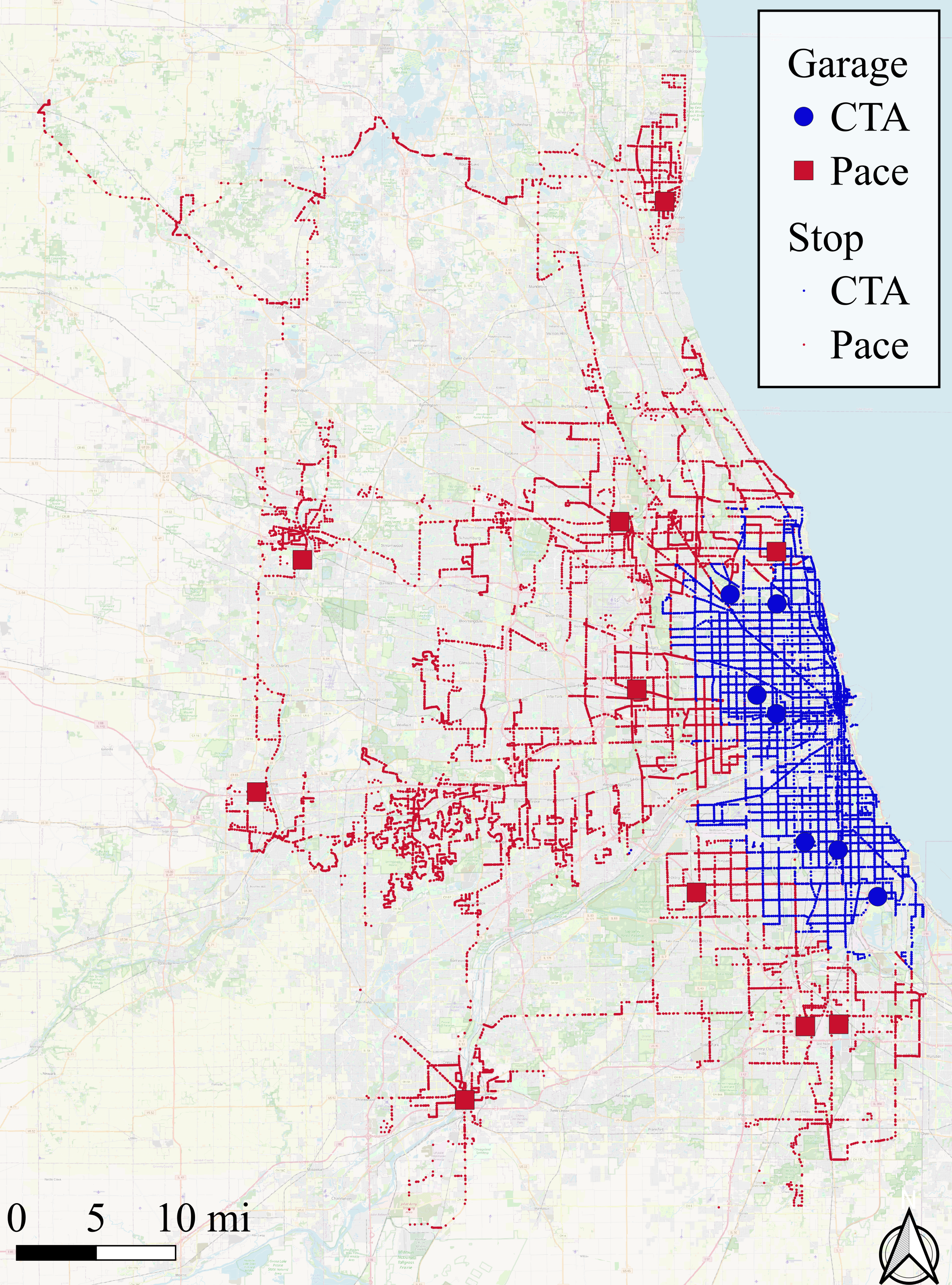} 
        \caption{Candidate charging stations}
        \label[fig]{stops_and_garages}
    \end{subfigure}

    \caption{CTA and Pace bus stops and bus garages in the Chicago metropolitan region}
    \label[fig]{study_region}
\end{figure}

\begin{figure}[!ht]
    \centering
    \begin{subfigure}[t]{0.48\textwidth}
        \centering
        \includegraphics[width=\linewidth]{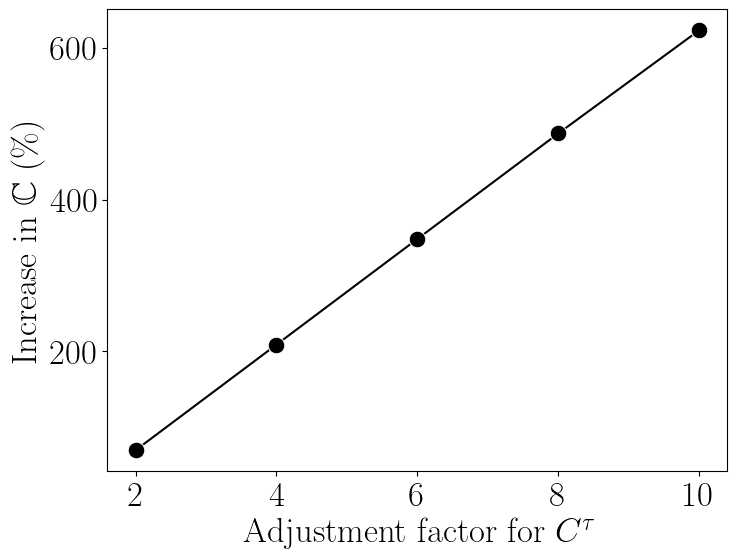} 
        \caption{Waiting time cost}
        \label[fig]{WCost_fig}
    \end{subfigure}
    \hfill
    \begin{subfigure}[t]{0.48\textwidth}
        \centering
        \includegraphics[width=\linewidth]{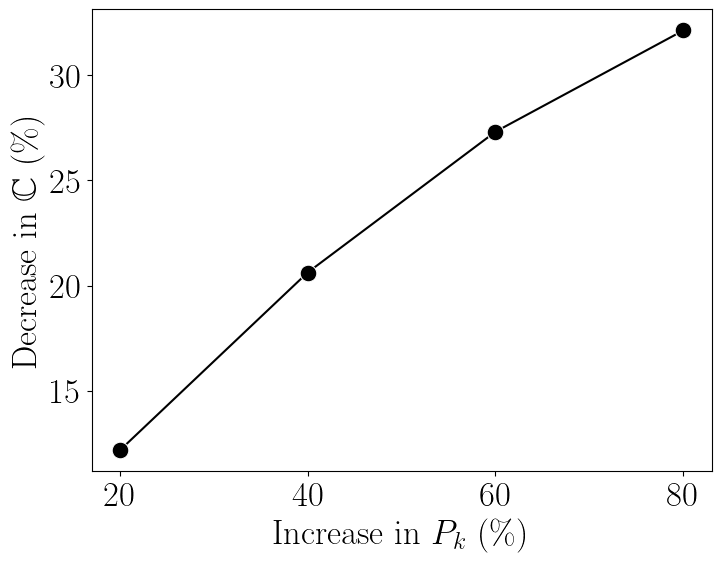} 
        \caption{Charger power}
        \label[fig]{ChPower_fig}
    \end{subfigure}

    \vspace{0.5cm}

    \begin{subfigure}[t]{0.48\textwidth}
        \centering
        \includegraphics[width=\linewidth]{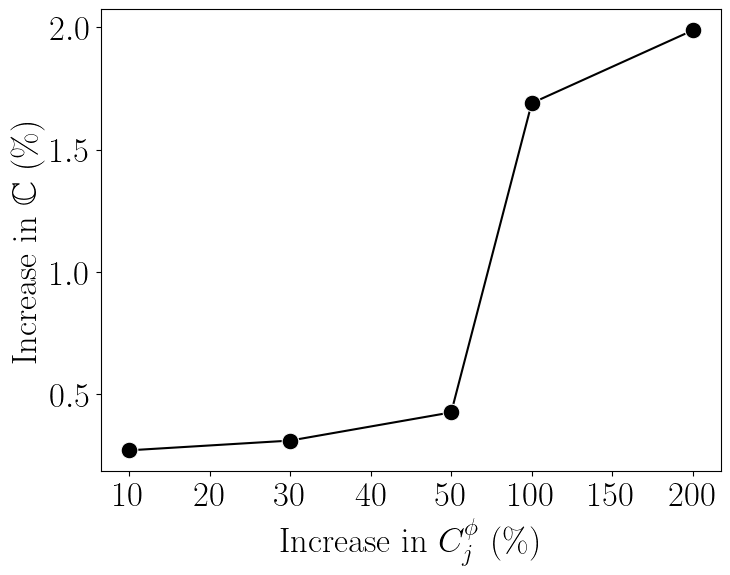} 
        \caption{Station cost}
        \label[fig]{SCost_fig}
    \end{subfigure}
    \hfill
    \begin{subfigure}[t]{0.48\textwidth}
        \centering
        \includegraphics[width=\linewidth]{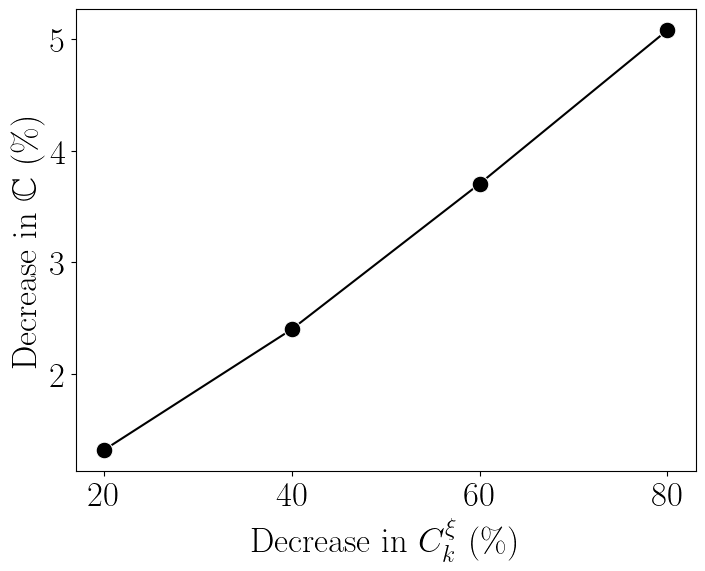} 
        \caption{Charger cost}
        \label[fig]{ChCost_fig}
    \end{subfigure}

    \caption{Average percentage change in objective value from baseline across scenarios.}
    \label[fig]{sensitivity_fig}
\end{figure}

\newpage
\section{Tables}\label[sup]{tables}

\begin{table}[!ht]
  \footnotesize
  \caption{A summary of the reviewed literature.}\label[tab]{lit_sum}
  \begin{center}
    \resizebox{\textwidth}{!}{\begin{tabular}{p{1.8cm}p{1.3cm} p{1.3cm} p{1.3cm} p{1.3cm} p{1.3cm} p{1.5cm} p{1.5cm} p{1.5cm}}
      \toprule
      Reference & Location & Allocation & Urban Transit & Charger Types & Waiting Time & Charging Type & Model & Case Study \\
      \midrule
      \cite{he2018optimal} & \cmark & \xmark & \xmark & \xmark & \xmark & terminals & bi-level & \xmark\\
      \cite{rogge2018electric} & \cmark & \cmark & \cmark & \xmark & \xmark & garages & MILP & Aachen \& Roskilde\\
      \cite{wei2018optimizing} & \cmark & \cmark & \cmark & \xmark & \xmark & terminals \& garages & MIP & Utah \\
      \cite{an2020battery} & \cmark & \xmark & \cmark & \xmark & \xmark & terminals & stochastic integer program & Melbourne \\
      \cite{USLU2021102645} & \cmark & \cmark & \xmark & \xmark & \cmark & terminals & MILP & Türkiye \\
      \cite{HSU2021103053} & \cmark & \xmark & \cmark & \xmark & \xmark & garages & augmented ILP & Taoyuan \\
      \cite{HaoHu2022} & \cmark & \xmark & \cmark & \xmark & \cmark & terminals & MILP & Sydney \\
      \cite{TZAMAKOS2023291} & \cmark & \cmark & \cmark & \cmark & \cmark & terminals & ILP & \xmark\\
      \cite{YiHe2023} & \xmark & \xmark & \cmark & \cmark & \xmark & terminals \& garages & MINLP & Salt Lake City\\
      \cite{DMcCabe2023} & \cmark & \cmark & \cmark & \xmark & \cmark & terminals \& garages & MILP & King County Metro \\
      \textbf{This paper} & \cmark & \cmark & \cmark & \cmark & \cmark & terminals \& garages & MINLP & Chicago \\
      \bottomrule
    \end{tabular}}
  \end{center}
  \emph{Note}: ILP: Integer Liner Programming, MILP: Mixed-Integer Linear Programming, MIP: Mixed Integer Programming
\end{table}

\begin{table}[!ht]
  \footnotesize
  \caption{Sets, parameters, and variables used in the model.}\label[tab]{notations}
  \begin{center}
    {\begin{tabular*}\textwidth{l p{12cm}}
      \toprule
      Set & Definition \\
      \midrule
      $\sI$ & set of charging demand occurrence points \\
      $\sI_j$ & subset of charging demand occurrence points from where BEBs can head to charging station $j \in \sJ$ \\
      $\sJ$ & set of candidate charging stations \\
      $\sJ_i$ & subset of charging stations that could be visited from charging demand occurrence point $i \in \sI$ \\ 
      $\sK$ & set of charger types \\
      \midrule
      Variable & Definition \\
      \midrule
      $s_{jk}$ & number of type $k$ chargers allocated to station $j \in \sJ$, $s_{jk} \in \mathbb{Z}_{\geq 0}$ \\
      $W_{jk}$ & expected waiting time for charger $k \in \sK$ at station $j \in \sJ$, $W_{jk} \in \mathbb{R}_{\geq 0}$ \\
      $x_{ijk}$ & $1$ if charger type $k \in \sK$ at station $j \in \sJ_i$ is visited before $i \in \sI$, $0$ otherwise \\
      $y_j$ & $1$ if station $j \in \sJ$ is an active station, $0$ otherwise \\
      \midrule
      Parameter & Definition \\
      \midrule
      $\epsilon$ & an infinitesimal number \\
      $\lambda_i$ & charging demand rate for point $i \in \sI$ \\
      $\mu_k$ & service rate of charger type $k \in \sK$ \\
      $C^\phi_j$ & fixed cost of charging station $j \in \sJ$ per time unit \\
      $C^\xi_k$ & cost of installing a type $k$ charger per time unit \\
      $C^\delta$ & cost of travel time per time unit \\
      $C^\tau$ & cost of waiting/charging time per time unit \\
      $T_{ij}^\delta$ & travel time from $i \in \sI$ to $j \in \sJ$ \\
      $P_k$ & charger power of charger type $k \in \sK$ \\
      $T_k^\rho$ & time needed for charger type $k \in \sK$ to recharge a BEB \\
      \bottomrule
    \end{tabular*}}
  \end{center}
\end{table}

\begin{table}[!ht]
  \footnotesize
  \caption{Parameter values in the baseline scenario.}\label[tab]{param_sum}
  \begin{center}
  \resizebox{\textwidth}{!}{%
    \begin{tabular}{llccr}
      \toprule
        \multirow{2.3}{*}{Parameter}  & \multirow{2.3}{*}{Definition} & \multicolumn{2}{c}{Value} & \multirow{2.3}{*}{Source}\\
        \cmidrule(lr){3-4}
        & & Slow & Fast & \\ 
        \midrule
        $P_k$ & charger power (kw) & 125  & 450 & \cite{authority2022charging, HaoHu2022} \\
        $C_k^\xi$ & charger cost (\$/$L^\xi$) & 55,500 & 200,000 & \cite{YiHe2023} \\
        $B$ & battery capacity (kWh) & \multicolumn{2}{c}{440} & \cite{authority2022charging} \\
        $L^\xi$ & charger lifetime (year) & \multicolumn{2}{c}{10} & \cite{an2020battery}\\
        $L^\phi$ & station lifetime (year) & \multicolumn{2}{c}{30} & \cite{HaoHu2022}\\
        $C_j^\phi$ & station cost (\$/$L^\phi$) & \multicolumn{2}{c}{208,000}  & \cite{YiHe2023}\\
        $C^\delta$ & travel cost (\$/min) & \multicolumn{2}{c}{2.67} & \cite{CTAAgencyProfile}\\
        $C^\tau$ & waiting/charging cost (\$/min) & \multicolumn{2}{c}{3.46} & \cite{AvgEnergyPrice}\\
      \bottomrule
    \end{tabular}
    }
  \end{center}
\end{table}

\newpage
\section{Algorithms}\label[sup]{algorithms}
\begin{algorithm}[!ht]
\caption{Smallest subset of stations for feasibility.}\label[algo]{minstations}
\SetArgSty{textnormal}
\scriptsize
\SetKwInOut{Input}{Input}
\SetKwInOut{Output}{Output}
\SetAlgoLined 
\SetKwFunction{FMinStations}{\textproc{MinStations}}
  \SetKwProg{Fn}{Function}{:}{\KwRet {$\sJ^{\min}$ \Comment{Smallest subset of stations for feasibility, and $\norm{\sJ^{\min}}$ is the minimum number of stations.}}}
  \Fn{\FMinStations{$\sI$, $\sI_j$, $\sJ$}}
  {
Initialize: $\sI^\psi \gets \sI$; $\sJ^\psi \gets \sJ$; $\sJ^{\min} \gets \emptyset$;

\While{$\sI^\psi \neq \emptyset$,}{
$j \gets \text{argmax}_{j\in\sJ^\psi}\norm{\sI_j}$; \Comment{Get station $j$ with the highest number of $i\in\sI_j$ ($j$ with min $C_j^\phi$ in tie).}

\If{$\sum_{k \in \sK} (\mu_k \overline{S}_{jk}) > \sum_{i \in \sI_j} \lambda_i$}{
\For{$i' \in \sI_j \cap \sI^\psi$}{
$\sI^\psi \gets \sI^\psi \setminus \{i'\}$; \Comment{Remove $i'$ from $\sI^\psi$.}
}
        $\sJ^{\min} \gets \sJ^{\min} \cup \{j\}$;\Comment{Add $j$ to $\sJ^{\min}$.}
        
        $\sJ^\psi \gets \sJ^\psi \setminus \{j\}$; \Comment{Remove $j$ from $\sJ^\psi$.}
}
}
}
\end{algorithm}

\begin{algorithm}[!ht]
\caption{Best number of type $k$ chargers allocated to station $j$ given $\sX$.}\label[algo]{s_given_x}
\SetArgSty{textnormal}
\scriptsize
\SetAlgoLined
\SetKwFunction{FBestChargers}{\textproc{BestChargers}}
  \SetKwProg{Fn}{Function}{:}{\KwRet {$s_{jk}^\prime$ \Comment{Best number of type $k$ chargers at $j$ given $\sX$.}}}
  \Fn{\FBestChargers{$\sX$}}
  {
  $s_{jk}^\prime \gets 0~\forall j\in\sJ,k\in\sK;$
  
\For{$(i,j,k)\in\sX$}{
$s_{jk}^{\min} \gets \sum_{i^\prime\in\sI_j|i=i^\prime}\frac{\lambda_{i^\prime}}{\mu_k}$;

\uIf{$s_{jk}^{\min} > \overline{S}_{jk}$}{
\Return$\emptyset$; \Comment{Solution $\sX$ is infeasible.}
}
\Else{
$s_{jk}^\prime \gets s_{jk}^{\min}$;

\While{$s_{jk}^\prime + 1 \leq \overline{S}_{jk}~\land~\sum_{i^\prime\in\sI_j|i=i^\prime}\lambda_{i^\prime} C^\delta \mathbb{W}(\cdot)^{+1} > C_k^\xi$
}{
\Comment{\parbox[t]{0.6\linewidth}{While incrementing $s_{jk}^\prime$ by 1 does not lead to infeasibility, and the decrease in $\mathbb{W}(\cdot)$ related costs with this increment is strictly greater than the cost of adding one more charger, do: increment $s_{jk}^\prime$ by 1.}}

$s_{jk}^\prime \gets s_{jk}^\prime + 1$
}
}
}
}
\end{algorithm}

\begin{algorithm}[!ht]
\caption{Assignment of each demand point $i$ to a charger type $k$ at an active charging station $j$ given the set of all active stations $\sJ^\alpha$.}\label[algo]{x_calc}
\SetArgSty{textnormal}
\scriptsize
\SetAlgoLined
\SetKwFunction{FDemandAssignment}{\textproc{DemandAssignment}}
  \SetKwProg{Fn}{Function}{:}{\KwRet {$\sX$ \Comment{The set of all tuples $(i,j,k)$ such that $y_j=1$ and $x_{ijk}=1$.}}}
  \Fn{\FDemandAssignment{$\sJ^\alpha$, $P$}}
  {
  $\sJ_i^\alpha \gets \sJ_i \cap \sJ^\alpha~\forall i\in\sI$; \Comment{\parbox[t]{0.6\linewidth}{$\sJ_i^\alpha$ is a subset of all active charging stations that could be visited from charging demand point $i \in \sI$.}}
  
\For{$i\in\sI$}{

$r \gets \text{\textproc{RandomUniform()}}$ \Comment{\textproc{RandomUniform()} generates a random number between 0 and 1.}

\uIf{$r < P$}{
$j \gets \text{random}({j\in\sJ^\alpha_i})$; \Comment{Assign the demand point $i \in \sI$ to a random open station.}
}
\Else{
$j \gets \text{argmin}_{j\in\sJ^\alpha_i}T^\delta_{ij}$; \Comment{Assign the demand point $i \in \sI$ to the closest open station.}
}
$k \gets \text{random}({k\in\sK})$; \Comment{Randomly select one of the charger types.}

$y_j \gets 1$; $x_{ijk} \gets 1$; $\sX \cup \{(i,j,k)\}$ \Comment{Add $(i,j,k)$ with $x_{ijk}=1$ to set $\sX$.}
}
}
\end{algorithm}

\begin{algorithm}[!ht]
\caption{Pseudocode for Simulated Annealing.}\label[algo]{SA_algo}
\SetArgSty{textnormal}
\scriptsize
\SetKwInOut{Input}{Input}
\SetKwInOut{Output}{Output}
\SetAlgoLined 
\SetKwFunction{FSimulatedAnnealing}{\textproc{SimulatedAnnealing}}
  \SetKwProg{Fn}{Function}{:}{\KwRet {$\sJ^{\beta}$, $\sX^\beta$, $s^\beta_{jk}$ \Comment{\parbox[t]{0.7\linewidth}{A set of active charging stations ($y_j=1$), assignment of each demand point to a charging station ($x_{ijk}=1$), and the number of each charger type at each facility ($s_{jk}$).}}}}
  \Fn{\FSimulatedAnnealing{$\sI$, $\sI_j$, $\sJ$, $P$, $T_0$, $L$, $C$}}
  {
Initialize: $T \gets T_0$, $l \gets 0$;

$\sJ^\theta \gets \text{\textproc{MinStations}}(\sI, \sI_j,\sJ)$; \Comment{\parbox[t]{0.45\linewidth}{Identify the subset of current active charging stations using \cref{minstations}.}}

$\sX^\theta \gets \text{\textproc{DemandAssignment}}(\sJ^\theta, P)$; \Comment{\parbox[t]{0.45\linewidth}{Calculate $\sX^\theta$ for the \textit{current stations} using \cref{x_calc}.}}

$s^\theta_{jk} \gets \text{\textproc{BestChargers}}(\sX^\theta) ~\forall j \in \sJ^\theta, k \in \sK$; \Comment{\parbox[t]{0.45\linewidth}{Calculate $s_{jk}$ for the \textit{current stations} using \cref{s_given_x}.}}

$\sJ^\beta \gets \sJ^\theta$; $\mathbb{C}^\beta \gets \mathbb{C}^\theta$; \Comment{\parbox[t]{0.45\linewidth}{Determine the objective value of the \textit{current stations} and set the \textit{best solution} to it.}}

\While{$l < L$}{

$l \gets l + 1$; \Comment{Increment the iteration count.}

$j \gets \text{\textproc{Random}}({j\in\sJ})$; \Comment{\textproc{Random}({$j\in\sJ$}) randomly selects one of the candidate stations.}

\uIf{$j \in \sJ^\theta$}{
    $\sJ^\nu \gets \sJ^\theta \setminus \{j\}$; \Comment{Create $\sJ^\nu$ by closing the active station $j \in \sJ^\theta$.}
}
\Else{
    $\sJ^\nu \gets \sJ^\theta \cup \{j\}$; \Comment{Create $\sJ^\nu$ by activating the closed station $j \in \sJ^\theta$.}
}

\While{$\sJ^\nu$ is not feasible,}{
$j \gets \text{\textproc{Random}}({j\in\sJ})$;

\uIf{$j \in \sJ^\nu$}{
    $\sJ^\nu \gets \sJ^\nu \setminus \{j\}$; \Comment{Remove $j$ from $\sJ^\nu$.}
}
\Else{
    $\sJ^\nu \gets \sJ^\nu \cup \{j\}$; \Comment{Add $j$ to $\sJ^\nu$.}
}
}
$\sX^\nu \gets \text{\textproc{DemandAssignment}}(\sJ^\nu, P)$; 

$s^\nu_{jk} \gets \text{\textproc{BestChargers}}(\sX^\nu) ~\forall j \in \sJ^\nu, k \in \sK$; 

\uIf{$\mathbb{C}^\nu < \mathbb{C}^\beta$}{
    $\sX^\beta \gets \sX^\nu$; $s^\beta_{jk} \gets s^\nu_{jk} ~\forall j \in \sJ^\nu, k \in \sK$,
    
    $\sJ^\beta \gets \sJ^\nu$; $\mathbb{C}^\beta \gets \mathbb{C}^\nu$; \Comment{Update the best solution with the new solution.}
    
    $\sJ^\theta \gets \sJ^\nu$; $\mathbb{C}^\theta \gets \mathbb{C}^\nu$; \Comment{Update the current solution with the new solution.}
}
\Else{
    $m = e^{\frac{\mathbb{C}^\theta - \mathbb{C}^\nu}{T}}$; \Comment{Calculate the metropolis acceptance criterion.}
    
    $r \gets \text{\textproc{RandomUniform()}}$; \Comment{\parbox[t]{0.55\linewidth}{\textproc{RandomUniform()} is a function that generates a random number between 0 and 1.}}

    \If{$r < m$}{
    $\sJ^\theta \gets \sJ^\nu$; $\mathbb{C}^\theta \gets \mathbb{C}^\nu$; \Comment{Update the current solution with the new solution.}
    }
}

$T \gets T \times \left(1 - \frac{C \times T_0}{L}\right)$; \Comment{\parbox[t]{0.6\linewidth}{Update temperature, and $C$ is a diminishing factor between 0 and 1 for decreasing temperatures.}}
}
}
\end{algorithm}

\begin{algorithm}[!ht]
\caption{$P$ subsets of stations with minimum sizes.}\label[algo]{covstations}
\SetArgSty{textnormal}
\scriptsize
\SetKwInOut{Input}{Input}
\SetKwInOut{Output}{Output}
\SetAlgoLined 
\SetKwFunction{FCoverSets}{\textproc{CoverSets}}
  \SetKwProg{Fn}{Function}{:}{\KwRet {$\sJ^{\sC}$ \Comment{\parbox[t]{0.75\linewidth}{Subset of stations ($\sJ^\alpha$) with $N$ number of subsets, and $\norm{\sJ^{\alpha}}$ is either the minimum or minimum+1 number of stations.}}}} 
  \Fn{\FCoverSets{$\sI$, $\sI_j$, $\sJ$, $N$}}
  {
Initialize $\sI^\psi \gets \sI$, $\sJ^\psi \gets \sJ$, $\sJ^{\sC} \gets \emptyset$, $S \gets \norm{\sJ}+1$;

\textproc{BackTrack}({$\sI^\psi$, $\emptyset$, $\sJ^\psi$, $\sJ^\sC$, $S$});
  }

\SetKwFunction{FBackTrack}{\textproc{Backtrack}}
  \SetKwProg{Fn}{Function}{:}{}
   \Fn{\FBackTrack{$\sI^\kappa$, $\sJ^\alpha$, $\sJ^\kappa$, $\sJ^\sC$, $S$}}{
\If{$\norm{\sJ^{\sC}} \geq N$}{
\Return \Comment{Population size is reached,}
}
\If{$\sI^\kappa = \emptyset$}{
    \uIf{$\norm{\sJ^\alpha} < S$}{
    $S \gets \norm{\sJ^\alpha}$; \Comment{Update minimum number of stations, $S$}

    $\sJ^{\sC} \gets \emptyset$; \Comment{Reset $\sJ^{\sC}$}

    $\sJ^{\sC} \gets \sJ^{\sC} \cup \sJ^\alpha$; \Comment{Add subset of active stations, ${\sJ^\alpha}$, to $\sJ^{\sC}$.}
    }
    \ElseIf{$\norm{\sJ^\alpha} = S~\lor~\norm{\sJ^\alpha} = S + 1$}{
    $\sJ^{\sC} \gets \sJ^{\sC} \cup \sJ^\alpha$; \Comment{Add subset of active stations, ${\sJ^\alpha}$, to $\sJ^{\sC}$.}
    }
    \Return
}
\If{$\norm{\sJ^\alpha} \geq S + 1$}{
    \Return \Comment{Remove the last added station to $\sJ^\alpha$ and check another station.}
}
\For{$j \in \sJ^\kappa$}{
\For{$i^\dagger \in \sI_j$}{
\If{$i^\dagger \in \sI^\kappa$}{
$\sI^{\kappa^\prime} \gets \sI^\kappa \setminus \{i^\dagger\}$; \Comment{Remove $i^\dagger\in\sI_j$ from new remaining demands $\sI^{\kappa^\prime}$ if $i^\dagger\in\sI^\kappa$.}
}
}
$\sJ^{\kappa^\prime} \gets \sJ^\kappa \setminus \{j\}$; \Comment{Remove station $j$ from new remaining stations $\sJ^{\kappa^\prime}$.}

$\sJ^{\alpha^\prime} \gets \sJ^\alpha \cup \{j\}$; \Comment{Insert station $j$ into new active stations $\sJ^{\alpha^\prime}$.}

\textproc{BackTrack}($\sI^{\kappa^\prime}$, $\sJ^{\alpha^\prime}$, $\sJ^{\kappa^\prime}$, $\sJ^\sC$, $S$); \Comment{Run the backtrack function with new inputs.}
}
}
\end{algorithm}

\begin{algorithm}[!ht]
\caption{Pseudocode for Genetic Algorithm.}\label[algo]{GA_algo}
\SetArgSty{textnormal}
\scriptsize
\SetKwInOut{Input}{Input}
\SetKwInOut{Output}{Output}
\SetAlgoLined 
\SetKwFunction{FGeneticAlgorithm}{\textproc{GeneticAlgorithm}}
  \SetKwProg{Fn}{Function}{:}{\KwRet {$\sJ^{\beta}$, $\sX^\beta$, $s^\beta_{jk}$ \Comment{\parbox[t]{0.7\linewidth}{A set of active charging stations ($y_j=1$), assignment of each demand point to a charging station ($x_{ijk}=1$), and the number of each charger type at each facility ($s_{jk}$).}}}}
  \Fn{\FGeneticAlgorithm{$\sI$, $\sI_j$, $\sJ$, $\sJ^\sC$,$P$, $L$}}
  {
\textbf{Initialize:} 

$l \gets 0$; \Comment{Initialize iteration count.}

$\sJ^\sC \gets \text{\textproc{CoverSets}}$($\sI$, $\sI_j$, $\sJ$, $N$); \Comment{Generate initial population \(\sJ^\sC\) using \cref{covstations}.}

$\sJ^\beta \gets \text{argmin}_{\sJ^\alpha_l\in\sJ^\sC}\mathbb{C}^\alpha_l$; \Comment{Identify the \textit{best solution} with the lowest objective value in \(\sJ^\sC\).}

\While{$l < L$}{

$l \gets l + 1$; \Comment{Increment the iteration count.}

\textbf{Selection:}

$\sJ^\sP \gets \text{\textproc{Random}}(\sJ^\sC, Y)$; \Comment{\textproc{Random}(\(\sJ^\sC, Y\)) selects a subset \(\sJ^\sP\) by choosing \(Y\%\) of \(\sJ^\sC\).}

$\sJ^{\pi1} \gets \text{argmin}_{\sJ^\alpha_l\in(\sJ^\sP)}\mathbb{C}^\alpha_l$; \Comment{Select the chromosome with the highest fitness as \textit{parent 1}.}

$\sJ^{\pi2} \gets \text{argmin}_{\sJ^\alpha_l\in{\sJ^\sP \setminus (\sJ^{\pi1})}}\mathbb{C}^\alpha_l$; \Comment{Select the second-best chromosome as \textit{parent 2}.}

\textbf{Crossover:}

$y_j^\nu=
\begin{cases} 
      y_j^{\pi1} & j\leq \frac{\norm{\sJ}}{2} \\
      y_j^{\pi2} & j > \frac{\norm{\sJ}}{2}
\end{cases}
\forall j \in \sJ$; \Comment{\parbox[t]{0.55\linewidth}{{Offspring} inherits genes from \textit{parent 1} for the first half of \(\sJ\), and from \textit{parent 2} for the second half.}}

$\sJ^\nu=\{j \mid y_j^\nu=1 ~\forall j\in\sJ\}$;

\textbf{Mutation:}

\While{$\sJ^\nu$ is not feasible,}{
$j \gets \text{\textproc{Random}}({j\in\sJ})$; \Comment{\textproc{Random}({$j\in\sJ$}) randomly selects one of the candidate stations.}

\If{$j \not\in \sJ^\nu$}{
    $\sJ^\nu \gets \sJ^\nu \cup \{j\}$; \Comment{Activate station \(j\) in \(\sJ^\nu\) if not already active.}
}
}

$\sX^\nu \gets \text{\textproc{DemandAssignment}}(\sJ^\nu, P)$; \Comment{Calculate $\sX^\nu$ for the {offspring} using \cref{x_calc}.}

\For{$(i,j,k) \in \sX^\nu$}{
\uIf{$j \leq \frac{\norm{\sJ}}{2} ~\land~ y_j^\nu = y_j^{\pi1} ~\land~ \sum_{k^\prime \in \sK}x_{ijk^\prime}^{\pi1}=1$}{
$k^{\pi1} \gets k^\prime \mid x_{ijk^\prime}^{\pi1}=1 \land k^\prime \in \sK$;

$\sX^\nu \setminus \{(i,j,k)\}$; $\sX^\nu \cup \{(i,j,k^{\pi1}\}$; \Comment{{Offspring} inherits charger types \(k\) from \textit{parent 1}.}
}
\ElseIf{$j > \frac{\norm{\sJ}}{2} ~\land~ y_j^\nu = y_j^{\pi2} ~\land~ \sum_{k^\prime \in \sK}x_{ijk^\prime}^{\pi2}=1$}{
$k^{\pi2} \gets k^\prime \mid x_{ijk^\prime}^{\pi2}=1 \land k^\prime \in \sK$;

$\sX^\nu \setminus \{(i,j,k)\}$; $\sX^\nu \cup \{(i,j,k^{\pi2}\}$; \Comment{{Offspring} inherits charger types \(k\) from \textit{parent 2}.}
}
}

$s^\nu_{jk} \gets \text{\textproc{BestChargers}}(\sX^\nu) ~\forall j \in \sJ^\nu, k \in \sK$; \Comment{\parbox[t]{0.4\linewidth}{Calculate $s_{jk}$ for the {offspring} using \cref{s_given_x}.}}

\textbf{Discard or add the offspring:}

$\sJ^\omega \gets \text{argmax}_{\sJ^\alpha_l\in\sJ^\sC}\mathbb{C}^\alpha_l$; \Comment{Identify the \textit{worst solution} with the highest objective value in \(\sJ^\sC\).}

\uIf{$\mathbb{C}^\nu < \mathbb{C}^\beta$}{
    $\sX^\beta \gets \sX^\nu$; $s^\beta_{jk} \gets s^\nu_{jk} ~\forall j \in \sJ^\nu, k \in \sK$,
    
    $\sJ^\beta \gets \sJ^\nu$; $\mathbb{C}^\beta \gets \mathbb{C}^\nu$; \Comment{Update the \textit{best solution}.}

    $\sJ^\sC \setminus \sJ^\omega$; $\sJ^\sC \cup \sJ^\nu$; \Comment{Replace the \textit{worst solution} $\sJ^\omega$ in \(\sJ^\sC\) with the {offspring} \(\sJ^\nu\).}
    }
\uElseIf{$\mathbb{C}^\nu > \mathbb{C}^\omega$}{
    $m = {\frac{\mathbb{C}^\omega}{\mathbb{C}^\nu}}$; \Comment{Calculate the acceptance criterion.}
    
    $r \gets \text{\textproc{RandomUniform()}}$ \Comment{\parbox[t]{0.55\linewidth}{\textproc{RandomUniform()} is a function that generates a random number between 0 and 1.}}

    \If{$r < m$}{
    $\sJ^\sC \setminus \sJ^\omega$; $\sJ^\sC \cup \sJ^\nu$; \Comment{Replace $\sJ^\omega$ in \(\sJ^\sC\) with $\sJ^\nu$.}
    }
}
\Else{
    $\sJ^\sC \setminus \sJ^\omega$; $\sJ^\sC \cup \sJ^\nu$; \Comment{Replace $\sJ^\omega$ in \(\sJ^\sC\) with $\sJ^\nu$.}
}
}
}
\end{algorithm}
\end{document}